\documentclass[letterpaper, paper,11pt]{AAS}		

\usepackage{AAS_packages}

\newcommand{\Poincare}{Poincar\'e }

\newcommand{\xk}{\ensuremath{x_k}}
\newcommand{\xkp}{\ensuremath{x_{k+1}}}
\newcommand{\yk}{\ensuremath{y_{k}}}
\newcommand{\ykp}{\ensuremath{y_{k+1}}}

\newcommand{\xdotk}{\ensuremath{\dot{x}_{k}}}
\newcommand{\ydotk}{\ensuremath{\dot{x}_{k}}}
\newcommand{\xdotkp}{\ensuremath{\dot{x}_{k+1}}}
\newcommand{\ydotkp}{\ensuremath{\dot{y}_{k+1}}}

\newcommand{\distonek}{\ensuremath{r_{1_k}}}
\newcommand{\distonekp}{\ensuremath{r_{1_{k+1}}}}
\newcommand{\disttwok}{\ensuremath{r_{2_{k}}}}
\newcommand{\disttwokp}{\ensuremath{r_{2_{k+1}}}}

\newcommand{\fonex}{\ensuremath{f_{1_x}}}
\newcommand{\ftwox}{\ensuremath{f_{2_x}}}
\newcommand{\fthreex}{\ensuremath{f_{3_x}}}
\newcommand{\ffourx}{\ensuremath{f_{4_x}}}

\newcommand{\foney}{\ensuremath{f_{1_y}}}
\newcommand{\ftwoy}{\ensuremath{f_{2_y}}}
\newcommand{\fthreey}{\ensuremath{f_{3_y}}}
\newcommand{\ffoury}{\ensuremath{f_{4_y}}}

\newcommand{\fonexd}{\ensuremath{f_{1_{\dot x}}}}
\newcommand{\ftwoxd}{\ensuremath{f_{2_{\dot x}}}}
\newcommand{\fthreexd}{\ensuremath{f_{3_{\dot x}}}}
\newcommand{\ffourxd}{\ensuremath{f_{4_{\dot x}}}}

\newcommand{\foneyd}{\ensuremath{f_{1_{\dot y}}}}
\newcommand{\ftwoyd}{\ensuremath{f_{2_{\dot y}}}}
\newcommand{\fthreeyd}{\ensuremath{f_{3_{\dot y}}}}
\newcommand{\ffouryd}{\ensuremath{f_{4_{\dot y}}}}

\PaperNumber{15-757}

\begin{document}

\title{Systematic Design of Optimal Low-Thrust Transfers for the Three-Body Problem}

\author{Shankar Kulumani\thanks{Doctoral Student, Mechanical and Aerospace Engineering, George Washington University, 800 22nd St NW, Washington, DC 20052, Email: \href{mailto:skulumani@gwu.edu}{skulumani@gwu.edu}.},  
Taeyoung Lee\thanks{Associate Professor, Mechanical and Aerospace Engineering, George Washington University, 800 22nd St NW, Washington, DC 20052, Tel: 202-994-8710, Email: \href{mailto:tylee@gwu.edu}{tylee@gwu.edu}.}
}

\maketitle{}

\begin{abstract}
A computational approach is developed for the design of continuous low thrust transfers in the planar circular restricted three-body problem.
The transfer design method of invariant manifolds is extended with the addition of continuous low thrust propulsion.
A reachable region is generated and it is used to determine transfer opportunities, analogous to the intersection of invariant manifolds.
The reachable set is developed on a lower dimensional \Poincare section and used to design transfer trajectories. 
This is solved numerically as a discrete optimal control problem using a variational integrator.
This provides for a geometrically exact and numerically efficient method for the motion in the three-body problem.
A numerical simulation is provided developing a transfer from a \( L_1 \) periodic orbit in  the Earth-Moon system to a target orbit about the Moon.
\end{abstract}

\section{Introduction}\label{sec:introduction}
Designing spacecraft trajectories is a classic and ongoing topic of research.
With the current fiscal constraints there is an increased focus on technologies with a critical impact on mass.
Optimal expenditure of onboard propellant is critical to allowing a mission to continue for a longer period of time or to enable the launch of a less massive spacecraft.
Electric propulsion systems offer a much greater specific impulse than chemical systems and are able to operate for extended periods of time.
However, these electric propulsion systems typically have much less thrust than their chemical counterparts and therefore must operate over longer durations in order to impart the desired momentum change.
Recent developments in miniature electric propulsion offer the potential for new research opportunities for small spacecraft~\cite{haque2013}.
With reduced development intervals and decreased launch costs, small satellites have gained increased attention as a cost effective means of scientific and technologic development. 
The merger of small satellites with miniature electric propulsion enables inexpensive and responsive missions requiring large changes in orbital energy or extended mission lifespan.
With the potential for more demanding missions, even greater importance is placed on the mission design to ensure that optimal trajectories satisfy mission requirements. 
In addition, non-Keplerian orbits and multi-body dynamics have been shown to allow for a much greater range of potential missions at a reduced energy cost~\cite{koon2000}.
Future space missions are increasing in complexity and will require new classes of orbits that are not possible via the traditional conic approach~\cite{ross2006,gomez2001}.
Optimally combining the dynamical structure of the three-body problem with low-thrust propulsion systems is vital for future mission success.

There has been extensive research focused on optimal control for spacecraft orbital transfers in the three-body problem~\cite{mingotti2011,grebow2011}.
Typically, the optimal control problem is solved via direct methods, which approximate the continuous problem as a parameter optimization problem.
The state and/or control trajectories are parameterized and solved in the form of a nonlinear optimization problem.
References~\citenum{mingotti2011} and~\citenum{grebow2011} use this direct approach in designing low-thrust transfers in the three-body problem.
Alternatively, indirect methods apply calculus of variations to derive the necessary conditions for optimality. 
This yields a lower dimensioned problem than the direct approach and algebraic conditions that, when satisfied, guarantee local optimality in contrast to direct methods which result in sub-optimal solutions.

The application of optimal control methods for orbital trajectory design is nontrivial.
The three-body system dynamics are nonlinear and exhibit chaotic behaviors. 
Small changes in initial conditions result in large variations of the resulting system trajectory. 
Therefore, any optimization routine is highly sensitive to the initial guess.
In addition, insight into the problem or intuition on the part of the designer is often required to determine initial conditions that will converge which achieve satisfactory results.
Efficient numerical implementation is dependent on correct initial conditions as well as accurate numerical integration.

Additionally, References~\citenum{mingotti2011} and~\citenum{grebow2011} implement the solutions using conventional Runge-Kutta integration techniques.
These techniques suffer from numerical instability and energy drift behaviors which make them ill-suited for long-term propagation.
These dissipative effects are even more detrimental with the addition of low-thrust propulsion to the dynamic equations of motion.
Conventional integration techniques fail to capture the physical laws and geometric properties of the dynamic system.
As a result, the long term effects of low-thrust on the spacecraft trajectory are not accurately captured. 

References~\citenum{koon2000} and~\citenum{ross2006} have illustrated the rich structure that exists in the three-body problem.
Within the three-body problem, a spacecraft's feasible region of motion is constrained by its energy, or Jacobi integral. 
It has been shown that there exist multi-dimensional tubes, or invariant manifolds, of constant energy trajectories that span the state space. 
Associated with periodic solutions of the three-body dynamics, these invariant manifolds allow for the spacecraft to traverse vast expanses of the state space with zero energy change. 
However, the results presented are highly case specific and difficult to generalize to arbitrary transfers.
Also, these results are based on control-free trajectories which rely on the underlying structure of the three-body system.
In addition, transfer orbits along an invariant manifold require a longer time of flight which may be undesirable for time critical missions.
The addition of low-thrust propulsion offers the potential of reduced transit times and the ability to depart from the free motion trajectory to allow for increased transfer opportunities. 

In order to address these issues, the authors propose an accurate and numerically stable method for long-term optimal orbital transfers.
This approach will avoid the instability and dissipative effects of conventional integration schemes.
In addition, the effects of low-thrust propulsion will be directly included in the integration method to ensure that extended duration optimal trajectories are computed correctly.
This improved method will enable the derivation of a systematic method of generating optimal transfer orbits between arbitrary states.
Indirect optimal control, based on the calculus of variations, will be implemented in order to generate optimal orbital transfers.
This will avoid the approximation issues inherent in the previous work, which utilized direct optimal control methods.
With this proposed method, the previous research on control-free trajectories will be generalized with the addition of low-thrust propulsion systems.

To achieve these objectives, computational geometric optimal control techniques are applied. 
The dynamics of the three-body system are derived from the discrete Lagrangian, which approximates the integral of the continuous time Lagrangian over a fixed discrete step.
Application of the discrete Euler-Lagrange equations, or the discrete Legendre transform, results in the discrete equations of motion.
This discrete update map, or variational integrator, shares the same geometric properties of the continuous time system and exhibits much better energy behavior than the traditional integration methods, especially over long time periods.
A discrete optimal control problem is formulated from the discrete equations of motion.
This approach, where explicit discretization occurs prior to optimization,  is in contrast to the typical method, where the equations of motion are implicitly discretized during the optimization procedure.
Formulating the problem in this manner results in more stable and accurate optimal solutions. 
In indirect methods the optimal control problem is expressed as a two-point boundary value problem.
Optimal solutions are generally sensitive to small variations in the initial multipliers.
As a result, the numerical stability of sensitivity derivatives is critical to accuracy and computational performance. 
The use of geometric integrators, which do not suffer the numerical dissipation of conventional integration methods, results in a more robust and efficient solution.

A discrete optimal control problem is formulated to determine the reachability set on a \Poincare section.
Given an initial condition and fixed time horizon, the reachable set is the set of states attainable, subject to the operational constraints of the spacecraft. 
Generation of this reachability set allows for the extension of the previous control-free missions in the three-body problem.
In addition, the generation of the reachable set allows for a more systematic method of determining initial conditions and eases the burden on the designer. 
The addition of low-thrust propulsion affords the ability to enlarge the reachable set as compared to the control-free case.
Maximization of the reachability set, on an appropriately chosen \Poincare section, allows for a greater space of potential transfer trajectories.
The use of the \Poincare section allows for design on a lower dimensional space and simplifies the design process.
Once an intersection is determined on the \Poincare section a transfer is calculated.

In short, the authors present a systematic method of generating optimal transfer orbits in the three-body problem.
Previous results in the design of optimal transfers have relied on suboptimal direct optimization methods and conventional integration techniques.
This paper provides a discrete optimal control formulation to generate the reachability set on a \Poincare section.
The use of a geometric integrators ensures numerical stability for long-duration orbit transfers.
A numerical example is presented demonstrating a transfer trajectory from \( L_1 \) to the vicinity of the Moon and compared to the control-free design methodology.
\section{System Model}\label{sec:pcrtbp}
The system model is based on the planar circular restricted three body problem~(PCRTBP).
The Earth is assumed to be the more massive primary, \( m_1 \), while the Moon is the second, smaller primary \( m_2\).
The equations of motion are developed in a rotating reference frame which allows for much greater insight into the structure of the dynamics.
The \( \hat{x} \) axis is directed along the vector from the Earth to the Moon.
The \( \hat{y} \) axis lies in the orbital plane and is orthogonal to \( \hat{x} \).
The rotating reference frame is centered at the system barycenter.
It is assumed that the \(\left( \hat{x}, \hat{y}\right)\) rotates with a constant angular velocity equal to the mean motion of the Moon.
Following convention, the system is also non-dimensionalized by the characteristic mass, length, and time~\cite{koon2000}.
As a result, the system can be characterized by a single mass parameter \( \mu \),
\begin{equation}
	\mu = \frac{m_2}{m_1+m_2} \, .
	\label{eq:mass_param}
\end{equation}
The larger primary, \(m_1\), is located at \( \left(  -\mu , 0 \right)\) and the smaller \( m_2\) is located at \( \left( 1-\mu , 0 \right)\).
In the rotating reference frame the Lagrangian is given by
\begin{equation}
	L = \frac{1}{2} \left( \left( \dot{x} -y \right)^2 + \left( \dot{y} + x \right)^2 \right) + \frac{1-\mu}{r_1} + \frac{\mu}{r_2}\, ,
	\label{eq:lagrangian}
\end{equation}
where the distance \(r_1\) and \(r_2\) of the spacecraft to each primary is defined in rotating coordinates as
\begin{align}
	r_1 &= \sqrt{\left( x + \mu\right)^2 + y^2}\, , \\
	r_2 &= \sqrt{\left( x - 1 + \mu\right)^2 + y^2}\, .
	\label{eq:distances}
\end{align}
Application of the Euler-Lagrange equations results in following equations of motion defined in the rotating reference frame
\begin{align}
	\ddot{x} - 2 \dot{y} + \frac{\partial}{\partial x} U &= u_x \, ,\nonumber \\
	\ddot{y} + 2 \dot{y} + \frac{\partial}{\partial y} U &= u_y \, ,
	\label{eq:cont_eom}
\end{align}
where the effective potential \( U\) is defined as
\begin{equation}
	U = \frac{1}{2} \left( x^2 + y^2\right) + \frac{1-\mu}{r_1} + \frac{\mu}{r_2}\, ,
	\label{eq:eff_pot}
\end{equation}
and the control inputs are defined as \( u_x\) and \(u_y\).
The state is defined as \( \bar{x} = \begin{bmatrix}\bar{r} &\bar{v} \end{bmatrix}\) with \(\bar{r} \in \R^{2\times1}\) and \(\bar{v} \in \R^{2\times1}\) representing the position and velocity with respect to the system barycenter, respectively.
The equations of motion may be rewritten in state space form as
\begin{equation}
	\left[\begin{array}{c} \dot{\bar{r}} \\ \dot{\bar{v}} \end{array} \right] = 
	\left[ \begin{array}{c} \bar{v} \\ A \bar{v} + \nabla U + \bar{u}(t) \end{array} \right] = f\left( t,x, u\right) \, ,
\end{equation}
where the matrix \( A \) and psuedo gravitational potential gradient \( \nabla U\) are
\begin{equation}\label{eq:A_mat}
	A = \left[ \begin{array}{ccc} 0 & 2 & 0 \\ -2 & 0 & 0 \\ 0 & 0 & 0 \end{array} \right] \, ,
\end{equation}
\begin{equation} \label{eq:grav_pot}
	\nabla U = \left[ \begin{array}{c} x - \frac{ \left(1 - \mu\right) \left(x + \mu\right)}{r_1^3} - \frac{\mu \left( x - 1 + \mu \right)}{r_2^3} \\
											y - \frac{ \left(1 - \mu\right) y}{r_1^3} - \frac{\mu y}{r_2^3} \\
											- \frac{ \left(1 - \mu\right) z}{r_1^3} - \frac{\mu z}{r_2^3}\end{array}\right]
					= \left[\begin{array}{c} U_x \\ U_y \\ U_z\end{array} \right] \, .
\end{equation}

\subsection{Jacobi Integral}
\begin{figure}[htbp]
	\centering
	\includegraphics[width=0.5\textwidth]{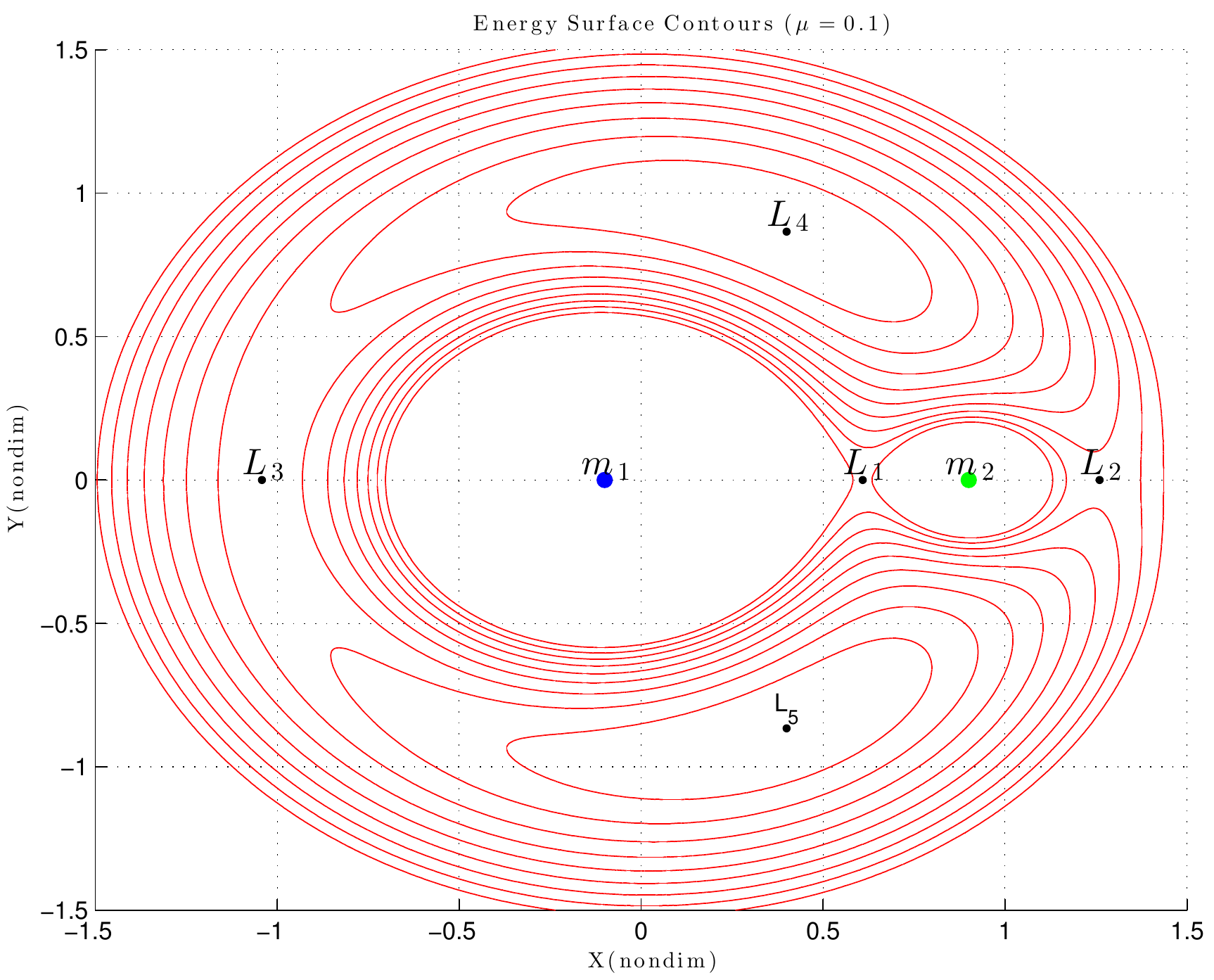}
	\caption{Contour Plot of Jacobi Integral}
	\label{fig:energy_contour}
\end{figure}

There exists a single integral, or constant of motion for the three-body problem~\cite{lanczos1970,szebehely1967}.
This energy constant is analogous to the total mechanical energy, however is a non-physical quantity arising from the problem formulation~\cite{szebehely1967}.
Also known as the Jacobi constant, it is defined as a function of the position and velocity in the rotating frame and given by
\begin{equation}
	E\left( \bar{r} , \bar{v} \right) = \frac{1}{2}\left( \dot{x}^2 + \dot{y}^2\right) - U\left(x,y \right) \, .
	\label{eq:jacobi}
\end{equation}
\Cref{eq:jacobi} divides the phase space into distinct regions based on the energy level of the spacecraft.
Fixing the Jacobi integral to a constant defines zero velocity curves, which are the locus of points where the kinetic energy, and hence velocity vanishes.
As seen in Figure~\ref{fig:energy_contour}, the phase space is divided into distinct realms based on the energy level.
In the vicinity of \( m_1\) or \(m_2\) there exists a potential well. 
As the energy level increases there are five critical points of the effective potential~\cref{eq:eff_pot} where the slope is zero.
Three collinear saddle points on the \(x\) axis and two equilateral points.
These equilibrium, or Lagrange points, are labeled \( L_i, i = 1, \hdots, 5 \) and are shown in Figure~\ref{fig:energy_contour}.
The Jacobi integral is a valuable invariant property of the three-body system that allows for greater insight into the motion of the spacecraft.

\section{Variational Integrator}\label{sec:discrete_var}
Geometric numerical integration deals with numerical integration methods which preserve the geometric properties of the flow of a differential equation, such as invaraint properties and symplecticity.
Variational integrators are constructed by discretizing Hamilton's principle rather than the continuous Euler-Lagrange equations~\cite{marsden2001}.
As a result, integrators developed in this manner have the desirable properties that they are symplectic and momentum preserving.
In addition, they exhibit improved energy behavior over long integration periods.
A short background on the variational principle for mechanical systems is presented. 
A discrete approximation of the action integral is presented and allows for construction of a variational integrator for the PCRTBP.

\subsection{Variational Principle}
Consider a continuous mechanical system described by the Lagrangian, \( L( q, \dot{q} ) \), for the generalized position, \( q\), and velocity, \( \dot{q} \).
In the standard approach of variational mechanics the action integral is formed by integrating the continuous Lagrangian along a path \( q(t) \) that the system follows from time \( t = 0 \) to \( t = T \)~\cite{greenwood1988}.
In the continuous time the action integral is defined as
\begin{align}\label{eq:action_integral}
	S = \int_{0}^T L\left( q, \dot{q}\right) \, dt \, .
\end{align}
Hamilton's principle states that the actual path followed by a holonomic system results in a stationary action integral with respect to path variations for fixed endpoints.
Taking the variation of~\cref{eq:action_integral} gives
\begin{subequations}\label{eq:var_principle}
\begin{align}
	\delta S &= \int_{0}^T \deriv{L}{q} \delta q + \deriv{L}{\dot{q}} \delta \dot{q} \, dt \\
		&= \int_{0}^T \deriv{L}{q} \delta q - \frac{d}{dt} \left( \deriv{L}{\dot{q}}\right) \delta q \, dt - \left. \left[ \deriv{L}{\dot{q}} \delta q\right] \right|_0^T \\
	&= \int_{0}^T \deriv{L}{q} - \frac{d}{dt} \left( \deriv{L}{\dot{q}}	\right) \, dt \, ,
\end{align}
\end{subequations} 
where we have used integration by parts and the conditions \( \delta q(0) = 0 \) and \( \delta q(T) = 0\).
For Hamilton's principle to be valid for all admissible variations \( \delta q \), the integrand of~\cref{eq:var_principle} must be zero for all \( t\), giving the continuous Euler-Lagrange Equations~\cite{lanczos1970}.
\begin{align}\label{eq:euler-lagrange}
	0 = \deriv{L}{q} - \frac{d}{dt} \left( \deriv{L}{\dot{q}} \right) \, .
\end{align}
Hamilton's equations are derivable through the use of the Legendre transformation which is a mapping \( \left( q, \dot{q},t\right) \rightarrow \left(q, p, t \right) \) where \( p_i\) is the generalized momenta.
\begin{align}\label{eq:legendre_transform}
	p_i = \deriv{L}{\dot{q}_i} \, .
\end{align}
In the continuous time case the Hamiltonian is defined as
\begin{align}\label{eq:hamiltonian}
	H &= \sum_{i = 1}^N p_i \dot{q}_i - L \left( q_i,\dot{q}_i, t \right) \, .
\end{align}
Applying~\cref{eq:legendre_transform} and taking the variation of~\cref{eq:hamiltonian} allows us to derive the equations of motion in Hamiltonian form
\begin{subequations}\label{eq:hamilton_eq}
\begin{align}
	\dot{q}_i &= \deriv{H}{p_i} \, ,\\
	\dot{p}_i &= - \deriv{H}{q_i} \, , \\
	\deriv{L}{t} &= -\deriv{H}{t} \, .
\end{align}
\end{subequations}
Both~\cref{eq:euler-lagrange,eq:hamilton_eq} result in equations of motion for the mechanical system and are equivalent via the Legendre transform.
\Cref{eq:euler-lagrange} results in \( n \) second order differential equations while~\cref{eq:hamilton_eq} results in \( 2n \) first order differential equations.
\subsection{Discrete Variational Mechanics}
A discrete analogue of Hamilton's principle and the action integral is formed.
Rather than taking a position, \( q \), and velocity, \( \dot{q} \), consider two positions \( q_0 \) and \( q_1 \) and a fixed time step \( h \in \R \).
The two positions are points on the curve \( q(t) \) such that \( q_0 \approx q(0) \) and \( q_1 \approx q(h) \).
A discrete time Lagrangian \( L_d( q_0, q_1) \) is formed which approximates the action integral between \( q_0 \) and \( q_1 \) as 
\begin{align}\label{eq:discrete_lagrangian}
	L_d\left( q_0 , q_1 \right) \approx \int_{0}^{h} L \left( q , \dot{q} \right) \, dt \, .
\end{align}
Since~\cref{eq:discrete_lagrangian} is calculated as a numerical integral, an appropriate quadrature rule is required.
There are multiple possible methods one can use to approximate the integral in~\cref{eq:discrete_lagrangian}.
An appropriate approximation rule is determined based on the ease of implementation and accuracy desired.
\begin{table}[htbp]
\caption{Selected Quadrature Rules\label{tab:quadrature}}
\begin{center}
\begin{tabular}{l|l}Rectangle & \( L_d(q_0,q_1) =L(q_0,\frac{q_1-q_0}{h}) h \)  \\ \hline
Midpoint & \( L_d(q_0,q_1) = L(\frac{q_0 + q_1}{2},\frac{q_1 - q_0}{h}) h \) \\ \hline
Trapezoidal & \( L_d(q_0, q_1) = \frac{1}{2} \left[ L(q_0, \frac{q_1 - q_0}{h} ) + L(q_1, \frac{q_1 - q_0 }{h} )\right] h \)
\end{tabular} 
\end{center}
\end{table}
\Cref{tab:quadrature} shows several possible approximation rules that are typically applied.
The rectangle rule is a first order accurate method and offers a straightforward implementation.
The midpoint and trapezoidal rules are both second order accurate methods. 
However, the midpoint rule results in an implicit form which adds further complexity to the equations of motion.
In this work, the trapezoidal approximation is applied to the PCRTBP.

Once an appropriate discrete Lagrangian is formed a discrete action sum is formed as the discrete analogue of~\cref{eq:action_integral}
\begin{align}\label{eq:action_sum}
	S_d = \sum_{k=0}^{N-1} L_d(q_k, q_{k+1}) \, .
\end{align}
Once again a discrete version of Hamilton's principle is applied to~\cref{eq:action_sum}.
Applying summation by parts yields
\begin{subequations}
\begin{align}\label{eq:dis_var_principle}
	\delta S_d &= \sum_{k=0}^{N-1} \deriv{L_d(q_k, q_{k+1})}{q_k} \delta q_k + \deriv{L_d(q_k, q_{k+1})}{q_{k+1}} \delta q_{k+1} \\
	&= \sum_{k=1}^{N-1} \bracket{ \deriv{L_d(q_k, q_{k+1})}{q_k} + \deriv{L_d(q_{k-1}, q_k)}{q_{k+1}}} \delta q_k \, .
\end{align}
\end{subequations}
For the discrete action sum to be stationary with respect to all admissible path variations, with fixed endpoints, the discrete Euler-Lagrange equations must be satisfied for \( k = 1, \cdots, N-1 \) resulting in
\begin{align}\label{eq:discrete_euler-lagrange}
	0 = \deriv{L_d(q_k, q_{k+1})}{q_k} + \deriv{L_d(q_{k-1}, q_k)}{q_{k+1}} \, .
\end{align}
A discrete version of the Legendre transformation, referred to as a discrete fiber derivative, results in the equivalent Hamiltonian form expression.
The discrete fiber derivative is given as 
\begin{subequations}\label{eq:discrete_legendre}
\begin{align}
	p_k &= \deriv{L_d(q_{k-1},q_k)}{q_{k+1}} = - \deriv{L_d(q_k, q_{k+1})}{q_k} \, ,\\
	p_{k+1} &= \deriv{L_d(q_k, q_{k+1})}{q_{k+1}} \, .
\end{align}
\end{subequations}
This yields a discrete Hamiltonian map \( (q_k, p_k) \to (q_{k+1}, p_{k+1}) \).
A more extensive development of variational integrators can be found in Reference~\citenum{marsden2001}.

\subsection{Discrete Equations of Motion}
The discrete equations of motion for the PCRTBP are derived by choosing an appropriate quadrature rule to discretize the Lagrangian in~\cref{eq:lagrangian}. 
In this work, the trapezoidal approximation is applied.
The trapezoid rule allows for an explicit second order accurate approximation.
The discrete Lagrangian is given by
\begin{align}\label{eq:discrete_lagrangian}
	L_d = &\frac{h}{2} \left( \frac{1}{2} \bracket{\left(  \frac{\xkp - \xk}{h} -\yk \right)^2 + \left( \frac{\ykp - \yk}{h} + \xk \right)^2} + \frac{1 - \mu}{r_{1_k}} + \frac{\mu}{r_{2_k}} \right. \nonumber \\ 
		& + \left. \frac{1}{2} \bracket{\left(  \frac{\xkp - \xk}{h} -\ykp \right)^2 + \left( \frac{\ykp - \yk}{h} + \xkp \right)^2} + \frac{1-\mu}{r_{1_{k+1}}} + \frac{\mu}{r_{2_{k+1}}}  \right) \, .
\end{align}
Applying a discrete version of the Lagrange-d'Alembert principle allows for inclusion of an external control force on the system~\cite{marsden2001}.
Using~\cref{eq:discrete_legendre,eq:discrete_lagrangian} and some manipulation, the equations of motion are given by
\begin{subequations}\label{eq:discrete_eoms}
\begin{align}
	\xkp &= \frac{1}{1+ h^2} \bracket{h \dot{x}_k + h^2 \dot{y}_k  + \xk \parenth{1+ \frac{3h^2}{2}} + \frac{h^3}{2} \yk - \frac{h^3}{2} U_{y_k} - \frac{h^2}{2} U_{x_k} } \label{eq:xkp} \, ,\\
	\ykp &= h \dot{y}_k + h \xk - h \xkp + \yk + \frac{h^2 \yk}{2} - \frac{h^2 }{2} U_{y_k} \label{eq:ykp} \, ,\\
	\dot{x}_{k+1} &= \dot{x}_k - 2 \yk + 2 \ykp + \frac{h}{2} \parenth{\xkp + \xk} - \frac{h}{2} U_{\xkp} - \frac{h}{2} U_{\xk} + h u_x \label{eq:xdotkp}\, ,\\
	\dot{y}_{k+1} &= \dot{y}_{k} + 2 \xk - 2 \xkp + \frac{h}{2} \parenth{\ykp + \yk} - \frac{h}{2} U_{\ykp} - \frac{h}{2} U_{\yk} + h u_y \label{eq:ydotkp} \, .
\end{align}
\end{subequations}
The discrete equations of motion are given in the Lagrangian form after applying the discrete fiber derivative from~\cref{eq:discrete_legendre} as \( p_{\xk} = \dot{x}_k - \yk \) and \( p_{\yk} = \dot{y}_k + \xk \).
The state is defined as \( \bar{ x}_k = \begin{bmatrix} \xk & \yk & \dot{x}_k & \dot{y}_k \end{bmatrix}^T\) and the control input is \( \bar{u} = \begin{bmatrix} u_x & u_y \end{bmatrix}^T \).
The discrete potential gradients are given by
\begin{subequations}\label{eq:discrete_potential_grad}
\begin{align}
	U_{\xk} &= \frac{\parenth{1 -\mu} \parenth{\xk + \mu}}{\distonek^3} + \frac{ \mu \parenth{\xk -1 + \mu}}{\disttwok^3} \label{eq:Uxk} \, ,\\
	U_{\yk} &= \frac{\parenth{1 -\mu} \yk}{\disttwok^3} + \frac{ \mu \yk}{\disttwok^3} \label{eq:Uyk} \, ,\\
	U_{\xkp} &= \frac{\parenth{1 -\mu} \parenth{\xkp + \mu}}{\distonekp^3} + \frac{ \mu \parenth{\xkp -1 + \mu}}{\disttwokp^3} \label{eq:Uxkp}\, ,\\
	U_{\ykp} &= \frac{\parenth{1 -\mu} \ykp}{\distonekp^3} + \frac{ \mu \ykp}{\disttwokp^3} \label{eq:Uykp} \, .
\end{align}	
\end{subequations}
The distances to each primary are defined as
\begin{subequations}\label{eq:discrete_distance}
\begin{align}
	\distonek &= \sqrt{\left( \xk + \mu\right)^2 + \yk^2} \label{eq:distonek}\, ,\\
	\disttwok &= \sqrt{\left( \xk - 1 + \mu\right)^2 + \yk^2} \label{eq:disttwok}\, ,\\
	\distonekp &= \sqrt{\left( \xkp + \mu\right)^2 + \ykp^2} \label{eq:distonekp}\, ,\\
	\disttwokp &= \sqrt{\left( \xkp - 1 + \mu\right)^2 + \ykp^2} \label{eq:disttwokp} \, .
\end{align}
\end{subequations}
Care must be taken during the implementation of~\cref{eq:discrete_eoms}.
As~\cref{eq:discrete_potential_grad,eq:discrete_distance} are defined at both step \( k \) and \( k+1 \) they must be evaluated at both time steps.
\Cref{eq:discrete_eoms} is implemented by first defining an initial state \( \bar{x}_k \) and control \( \bar{u}_k \).
The distances and gravitational potential at step \( k \) are evaluated from~\cref{eq:distonek,eq:disttwok,eq:Uxk,eq:Uyk}.
The discrete update steps in~\cref{eq:xkp,eq:ykp} are evaluated to generate \( \xkp \) and \( \ykp\).
Next, the distances and gravitational potential at step \( k+1 \) are evaluated from~\cref{eq:distonekp,eq:disttwokp,eq:Uxkp,eq:Uykp}. 
Finally, the update steps in~\cref{eq:xdotkp,eq:ydotkp} are evaluated.
This results in the complete discrete update map \( \bar{x}_k \to \bar{x}_{k+1} \) given \( \bar{u}_k \).

\begin{figure} 
	\centering 
	\begin{subfigure}[h]{0.3\textwidth} 
		\includegraphics[width=\textwidth]{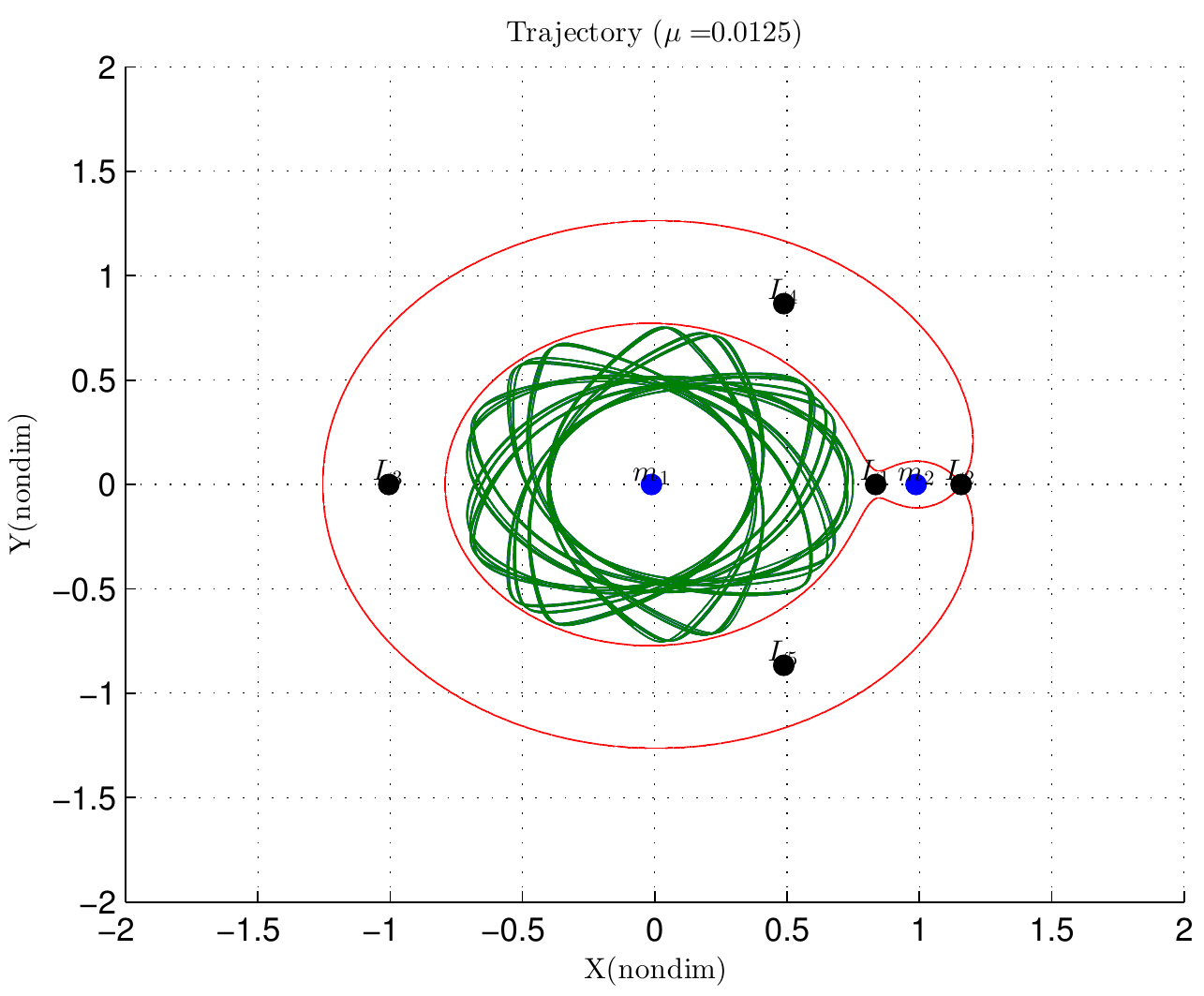} 
		\caption{Trajectory} \label{fig:compare_trajectory} 
	\end{subfigure}~ 
	\begin{subfigure}[htbp]{0.3\textwidth} 
		\includegraphics[width=\textwidth]{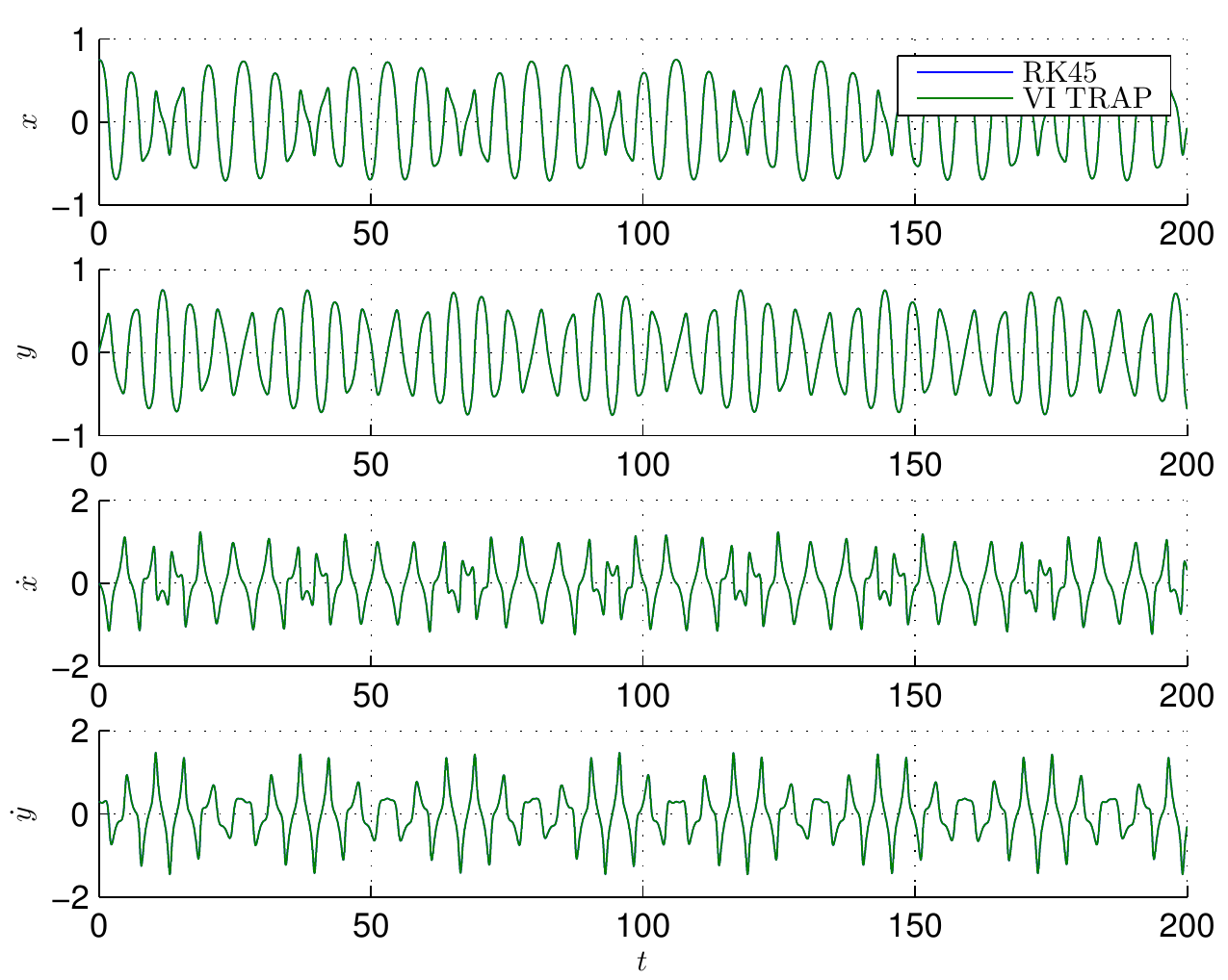} 
		\caption{State \( \parenth{x , y, \dot{x}, \dot{y}}\)} \label{fig:compare_components} 
	\end{subfigure} ~ 
	\begin{subfigure}[htbp]{0.3\textwidth} 
		\includegraphics[width=\textwidth]{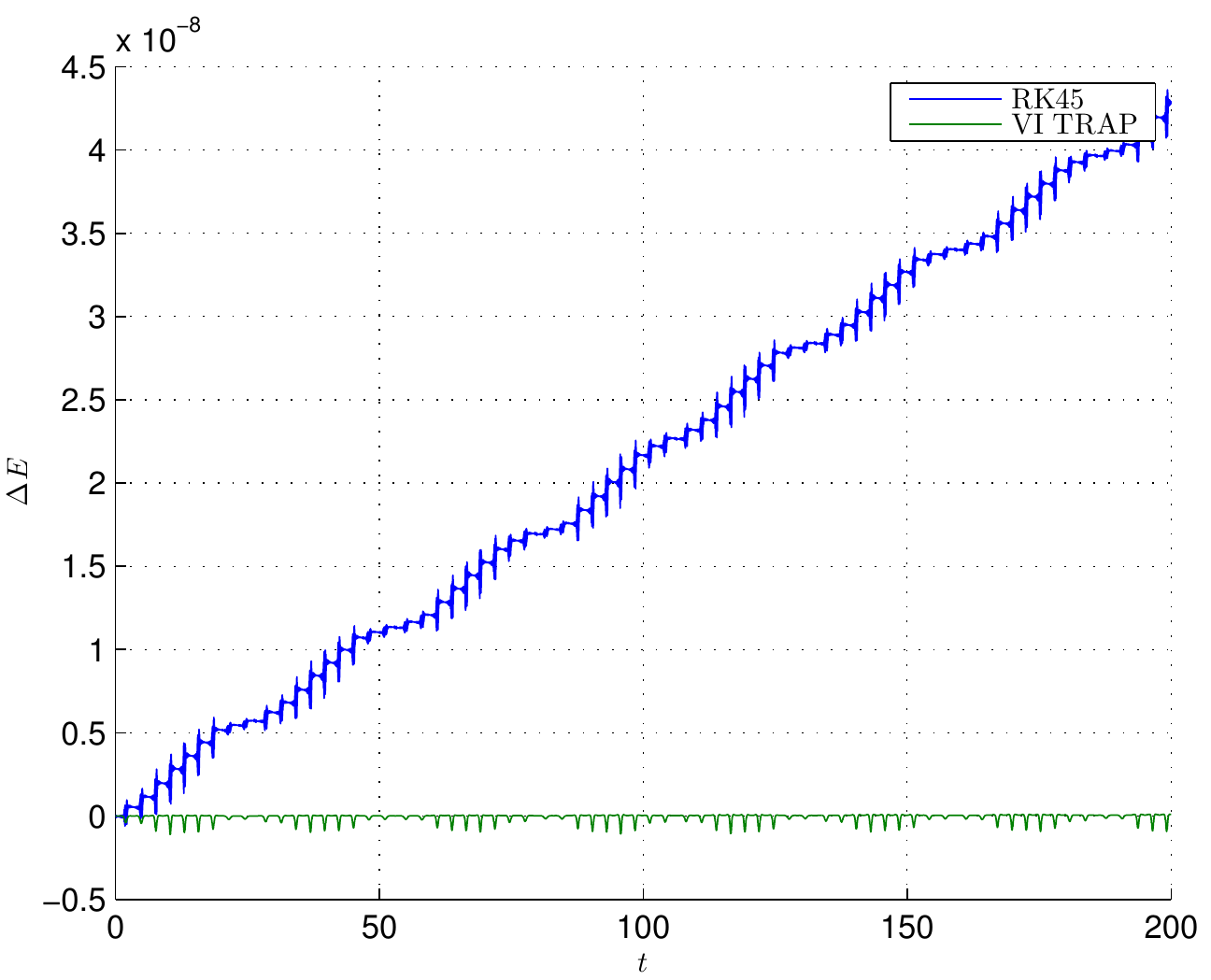} 
		\caption{Jacobi Integral} \label{fig:compare_energy} 
	\end{subfigure} 
	\caption{Integrator Comparison}
	\label{fig:integrator_compare} 
\end{figure}
A simulation comparing the variational integrator to a conventional Runge-Kutta method is given in~\cref{fig:integrator_compare}.
A particle is simulated from an initial condition of \( \bar{x}_0 = \begin{bmatrix} 0.75 & 0 & 0 & 0.2883\end{bmatrix}^T \) for \( t_f = 200 \approx 15\) years in the Earth-Moon system.
The variational integrator uses a step size of \SI{47.22}{\second} while the Runge-Kutta method uses a variable step size implemented via ODE45 in Matlab.
\Cref{fig:compare_trajectory} shows the trajectory of the spacecraft in the rotating reference frame for both integration schemes.
Both integration schemes result in trajectories that are initially nearly identical.
The discrete equations of motion are an accurate approximation for the continuous dynamics as they closely match the solution of ODE45 over the initial portion of the simulation.
However, as time progresses the trajectories begin to diverge due to the differences in system energy.
\Cref{fig:compare_energy} shows the evolution of the Jacobi integral.
The variational integrator exhibits a bounded behavior about the initial energy with a mean variation of \num{4.2522e-20}.
However, the conventional Runge-Kutta method demonstrates a clear energy drift of \num{4.2814e-8}. 
Over long simulation horizons or with the addition of small control inputs this poor energy behavior limits the applicability of conventional techniques.

\section{Invariant Manifolds} \label{sec:invariant_manifold}
\begin{figure}
     \centering
        \begin{subfigure}[b]{0.25\textwidth}
                \includegraphics[width=\columnwidth]{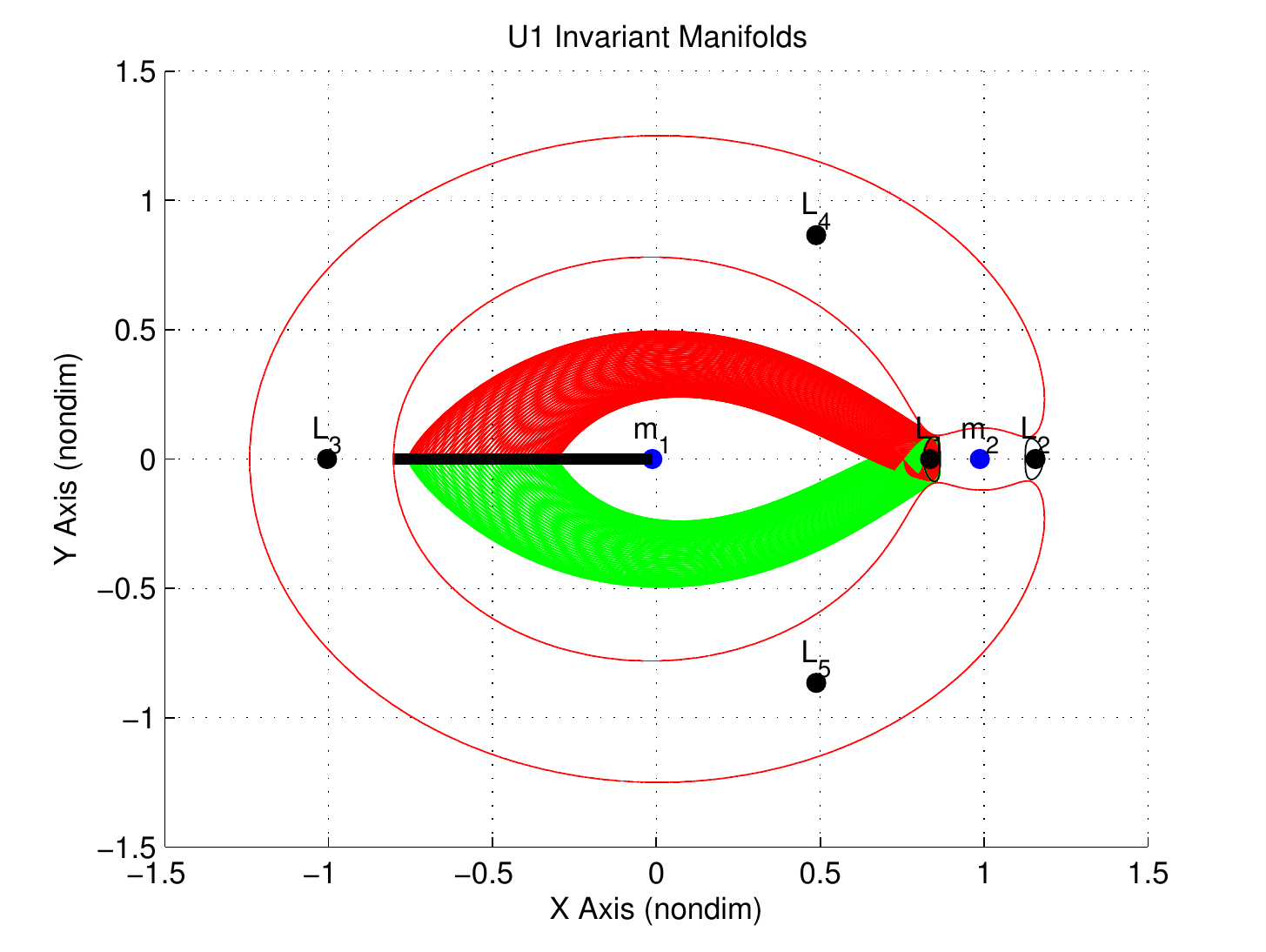}
                \caption{U1 Manifolds}
                \label{fig:u1_manifolds}
        \end{subfigure}%
        ~
        \begin{subfigure}[b]{0.25\textwidth}
                \includegraphics[width=\columnwidth]{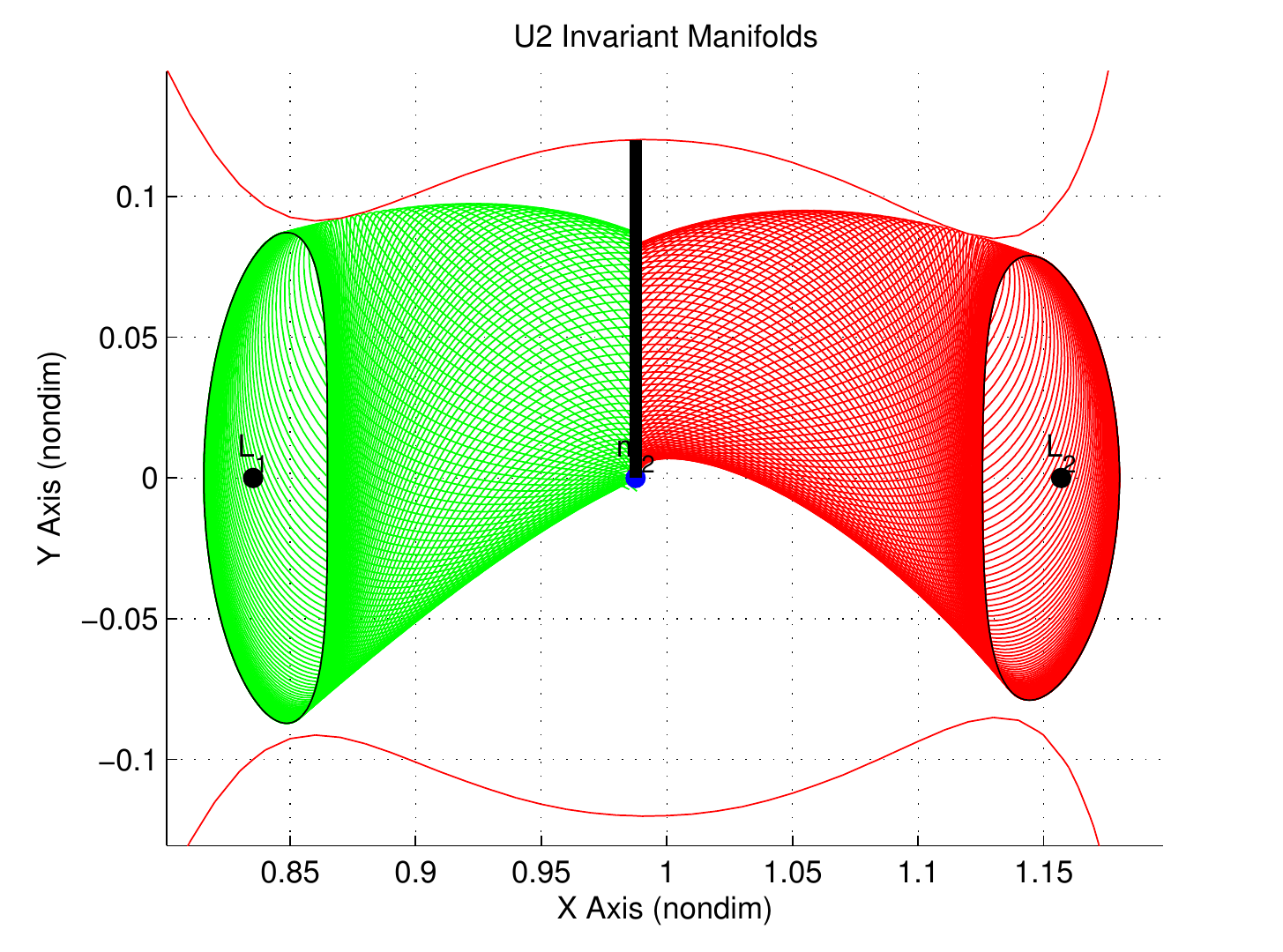}
                \caption{U2 Manifolds}
                \label{fig:u2_manifolds}
        \end{subfigure}~
        \begin{subfigure}[b]{0.25\textwidth}
                \includegraphics[width=\columnwidth]{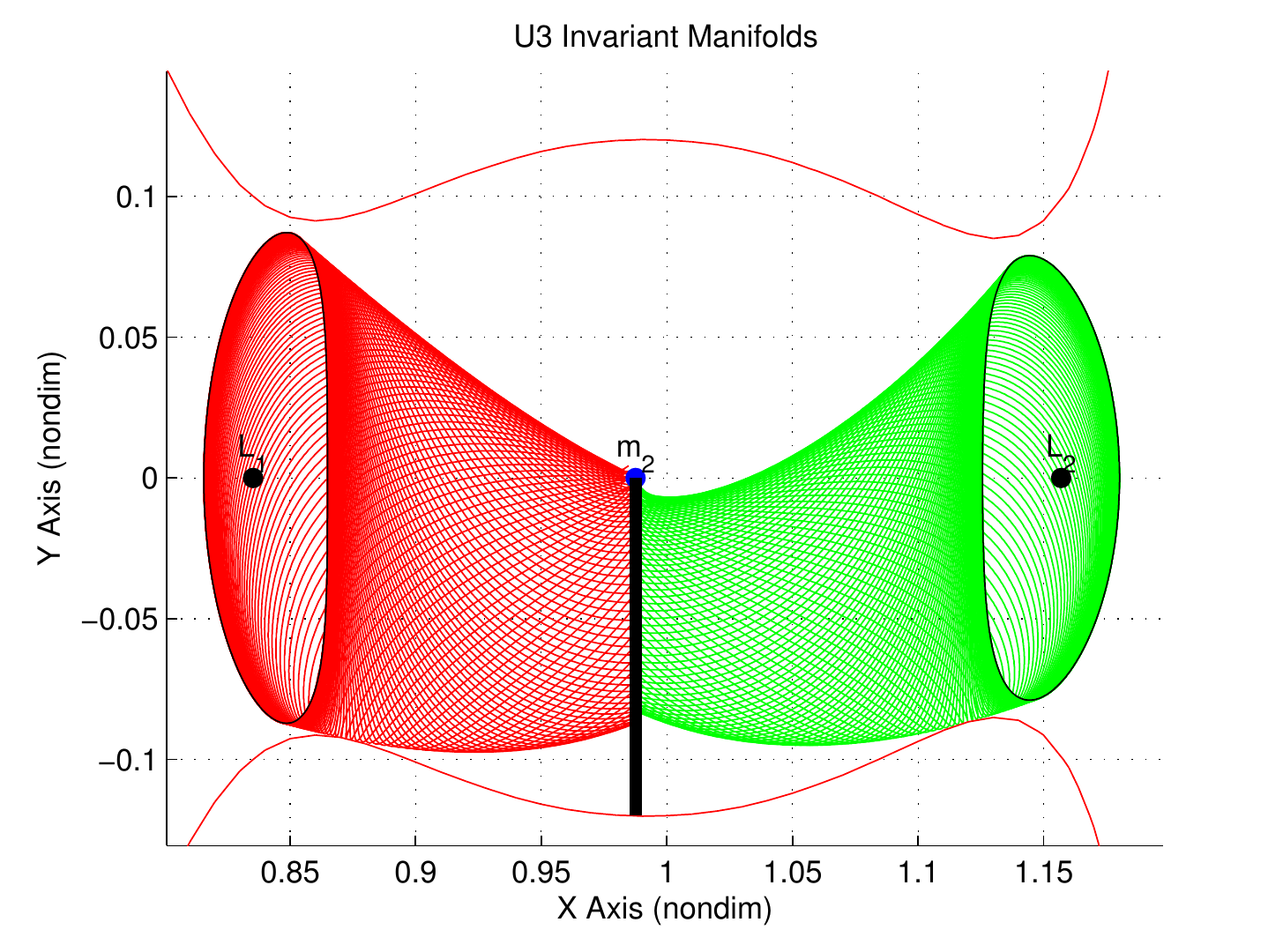}
                \caption{U3 Manifolds}
                \label{fig:u3_manifolds}
        \end{subfigure}%
        ~
        \begin{subfigure}[b]{0.25\textwidth}
                \includegraphics[width=\columnwidth]{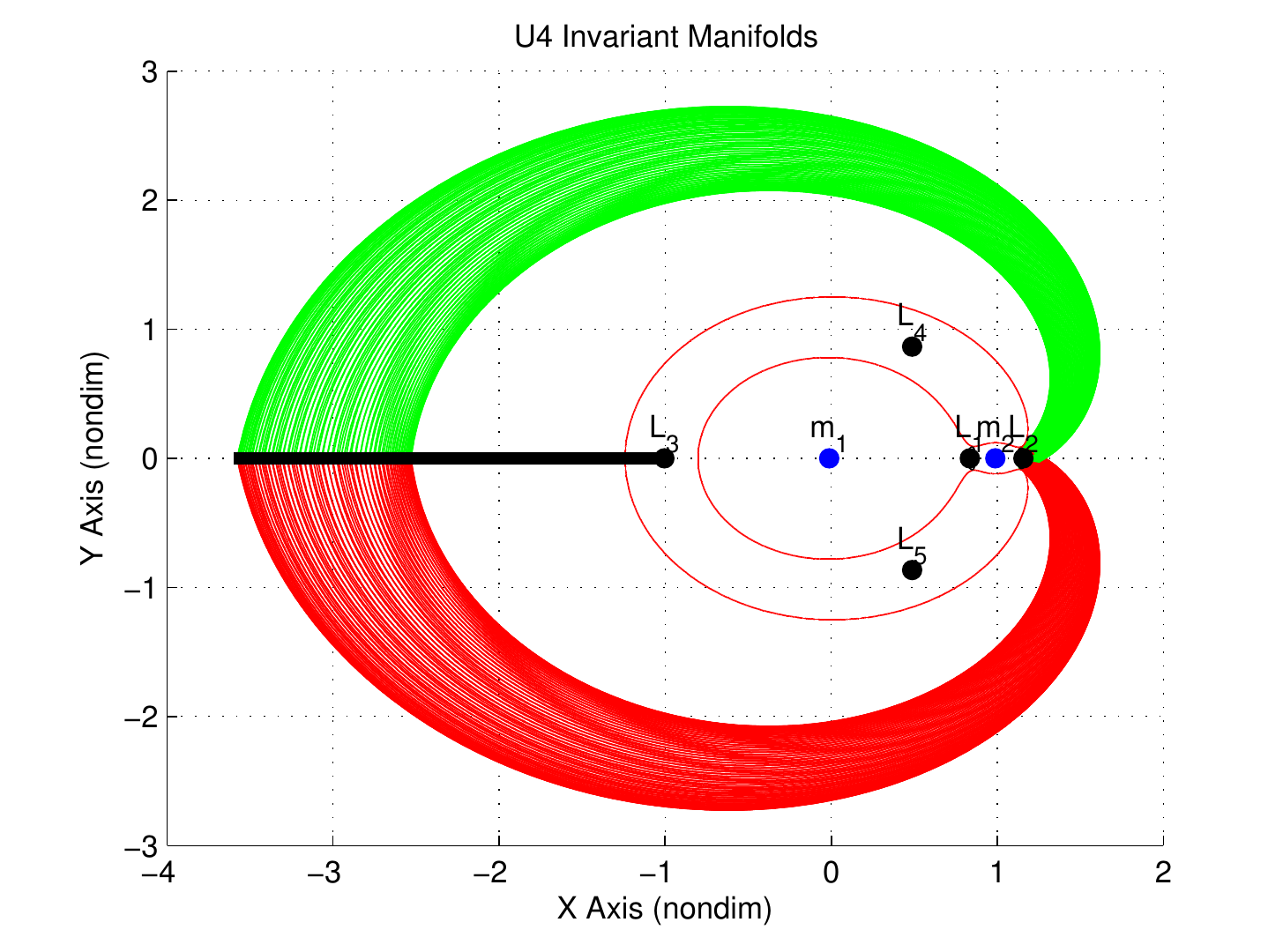}
                \caption{U4 Manifolds}
                \label{fig:u4_manifolds}
        \end{subfigure}
        \caption{Invariant Manifolds for Planar Earth-Moon three-body system}
	\label{fig:invariant_manifolds}
\end{figure}

There exists a rich structure of tubes, or invariant manifolds, in the three-body problem~\cite{koon2000,conley1968}.
The manifold structure associated with periodic orbits about the \( L_1 \) and \( L_2 \) Lagrange points are critical to the understanding of the motion of spacecraft as well as comets/asteroids.
In addition, the stable and unstable manifolds serve as the boundaries of the phase space region that allow for the transport between realms in a single three-body system or between multiple three-body systems.
These invariant manifolds only exist as a result of the dynamic formulation of the three-body problem in a rotating reference frame. 
In this way, it is possible to design trajectories that intersect the invariant manifolds and connect widely separated regions of the state space. 
\Cref{fig:invariant_manifolds} shows some of the tubes, projected from the phase space onto the configuration space, associated with \( L_1 \) and \( L_2 \) periodic orbits in the Earth-Moon system. 
This technique has been heavily investigated in Reference~\citenum{koon2000}, applied to operational missions in Reference~\citenum{koon1999}, and applied to potential multi-moon orbiter missions in Reference~\citenum{tanaka2011}.
In addition, much recent work has focused on potential return missions to the Moon~\cite{zanzottera2012,campagnola2012,mingotti2011,ozimek2010a,mingotti2009}.
\begin{figure}
     \centering
        \begin{subfigure}[b]{0.25\textwidth}
                \includegraphics[width=\columnwidth]{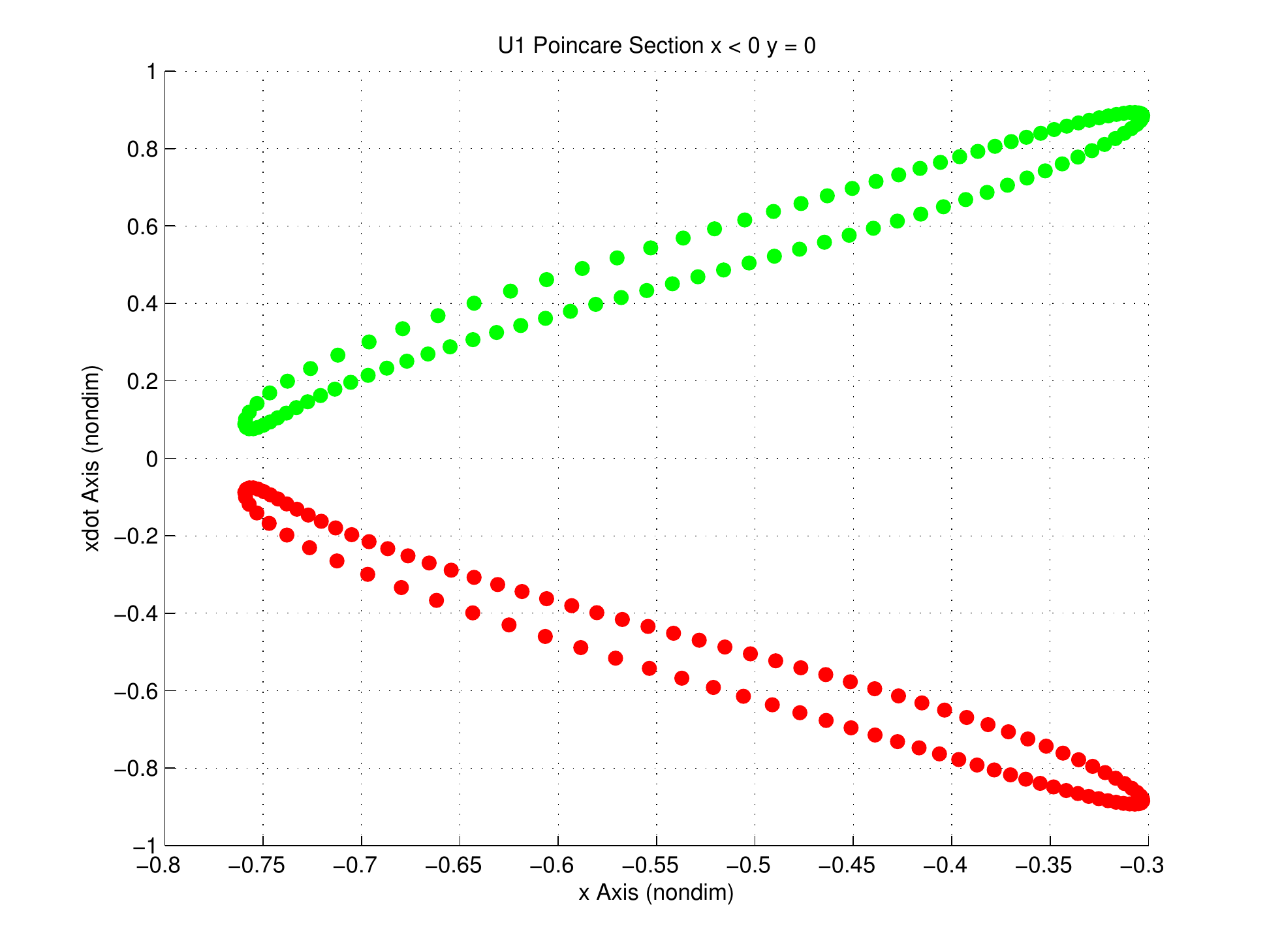}
                \caption{U1 Section}
                \label{fig:u1_poincare}
        \end{subfigure}%
        ~
        \begin{subfigure}[b]{0.25\textwidth}
                \includegraphics[width=\columnwidth]{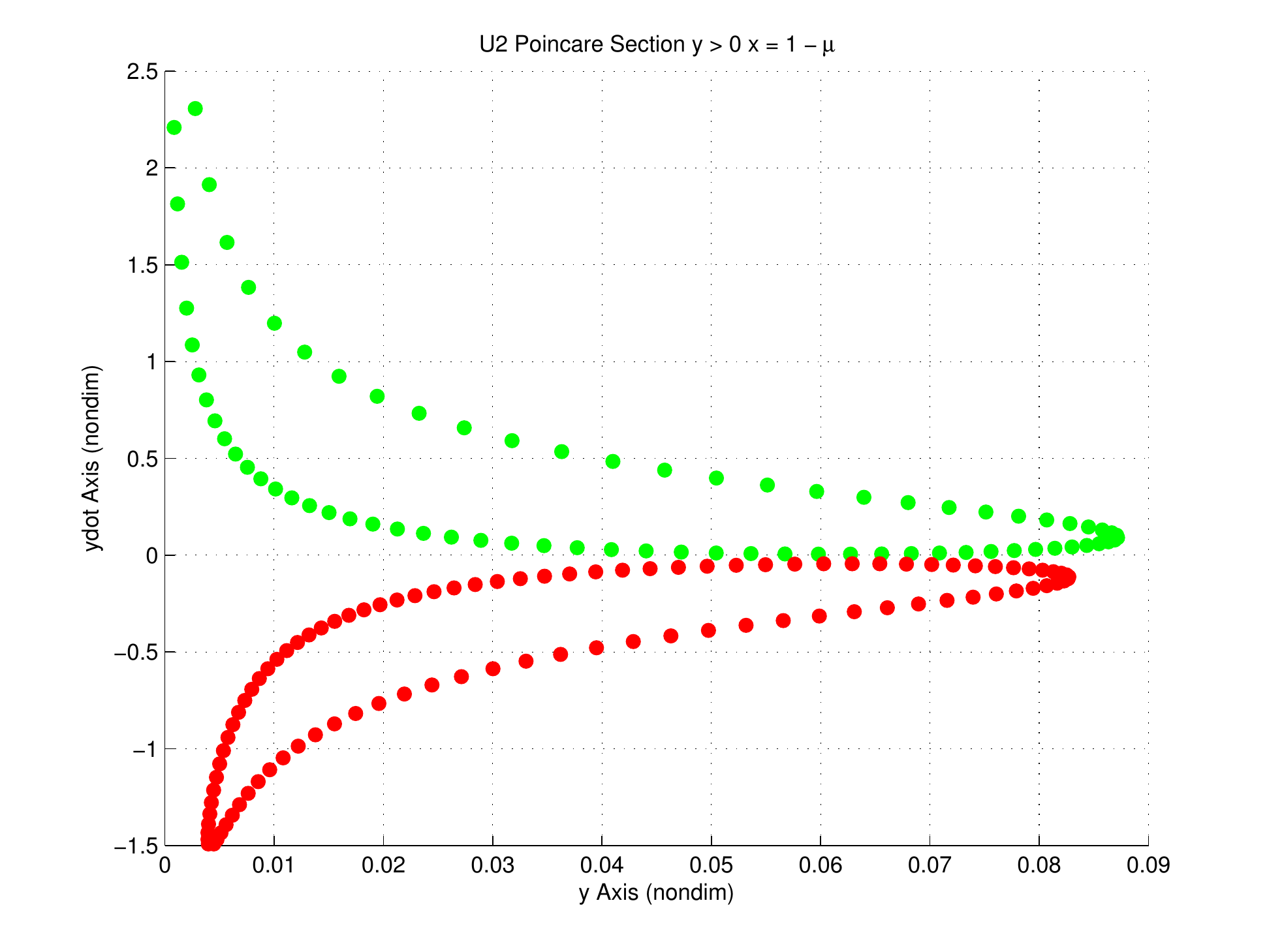}
                \caption{U2 Section}
                \label{fig:u2_poincare}
        \end{subfigure}~
        \begin{subfigure}[b]{0.25\textwidth}
                \includegraphics[width=\columnwidth]{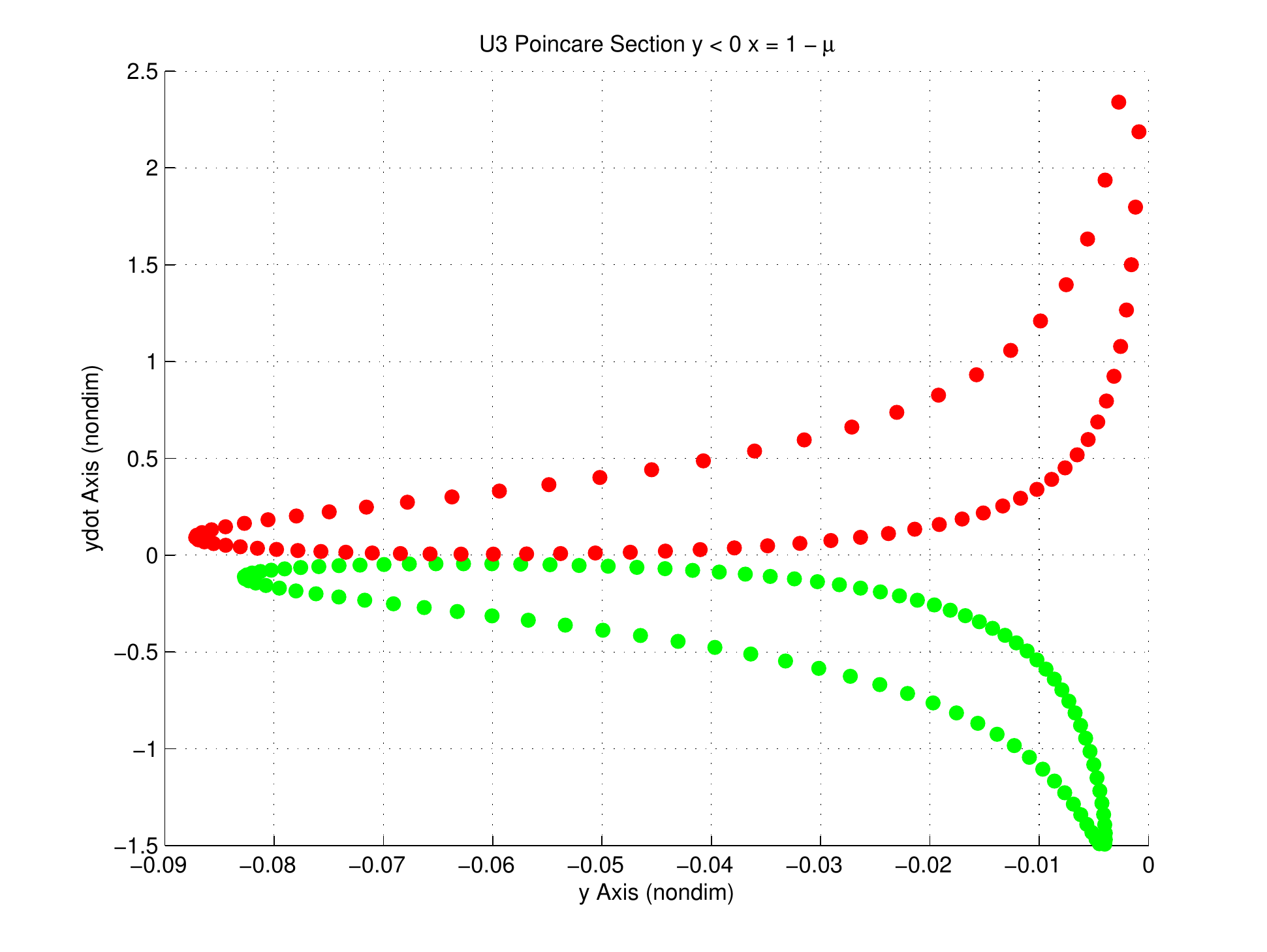}
                \caption{U3 Section}
                \label{fig:u3_poincare}
        \end{subfigure}%
        ~
        \begin{subfigure}[b]{0.25\textwidth}
                \includegraphics[width=\columnwidth]{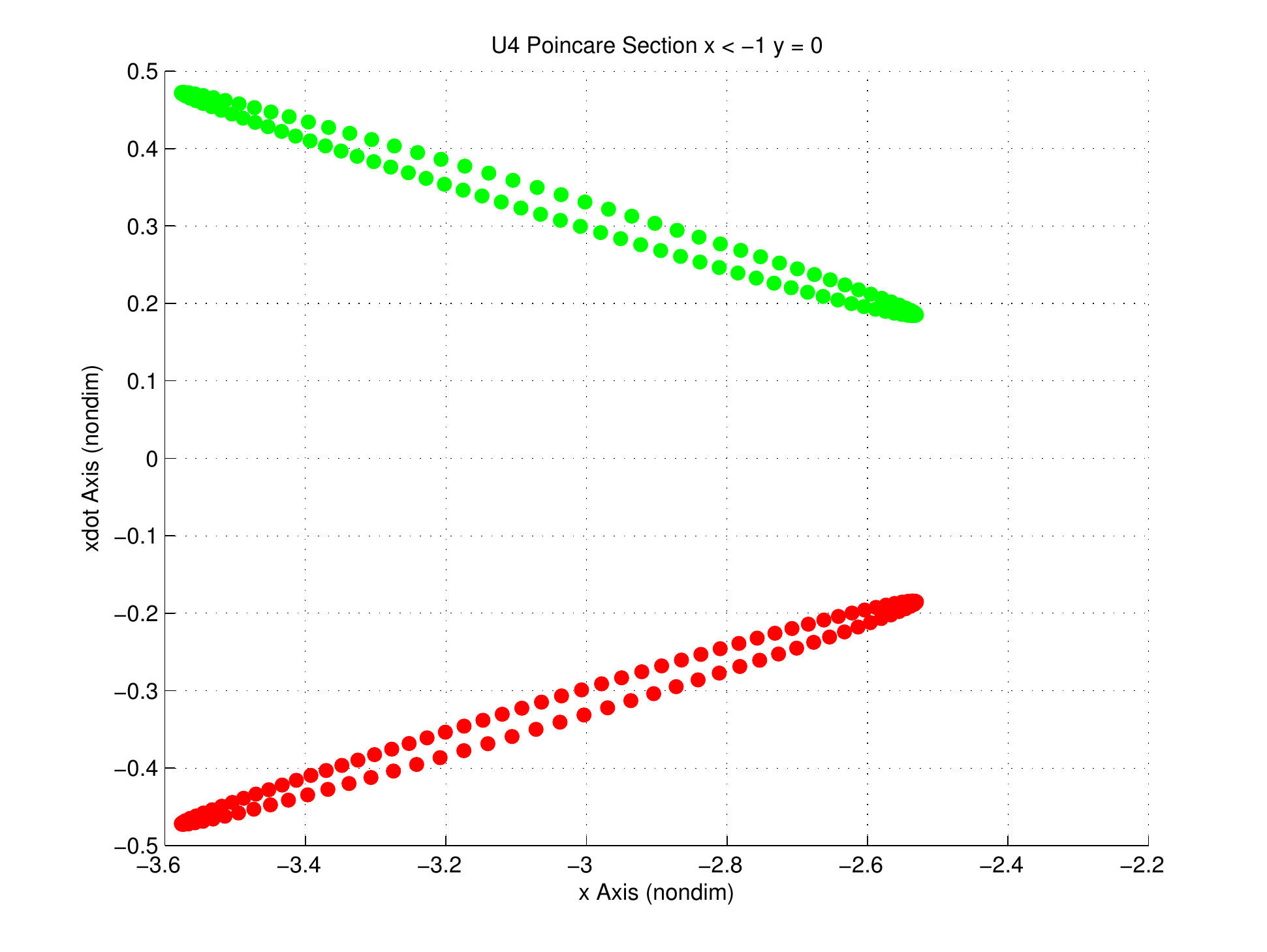}
                \caption{U4 Section}
                \label{fig:u4_poincare}
        \end{subfigure}
        \caption{\Poincare Sections for Planar Earth-Moon three-body system}
	\label{fig:poincare_sections}
\end{figure}

Another useful technique in the analysis of the free motion of dynamical systems is that of the \Poincare section.
The \Poincare section is the intersection of a periodic orbit of a dynamical system with a lower dimensional sub-space transverse to the flow.
This allows for greater insight into the stability and dynamics of periodic solutions of dynamic system. 
Points are drawn as the periodic solution intersects the \Poincare section.
Great insight and structure into the dynamical system is typically available through use of \Poincare section.
For example, \cref{fig:poincare_sections} shows the intersections of the invariant manifolds of~\cref{fig:invariant_manifolds} with the \Poincare sections defined by the black line segments.
The \Poincare section and the Jacobi integral reduce the state space from four to two dimensions for the planar three-body problem~\cite{koon2001}.
This greatly eases the design process and allows for a simple planar visualization of the intersection of the invariant manifolds.
Intersecting states are easy to determine and allow for transfers between invariant manifolds.

The intersection of the invariant manifolds on a \Poincare section is central to determining the desired control-free trajectory.
However, the results previously developed are highly case specific and difficult to generalize to arbitrary transfers.
Also, these results are based on control-free trajectories which rely on the underlying structure of the three-body system.
In addition, transfer orbits along an invariant manifold require large time of flights which may be undesirable for time critical missions.
The addition of low-thrust propulsion offers the potential of reduced transit times and the ability to depart from the free motion trajectory to allow for increased transfer opportunities. 
In this work, an optimal control problem is formulated to generate a reachable set of the spacecraft.
The reachable set is computed on an appropriate \Poincare section and is used to design an appropriate transfer trajectory.
\section{Optimal Control Formulation for Reachability Set}\label{sec:optimal_control}
In this section, an optimal control formulation is presented to determine and design transfers within the three-body problem.
The application of variational integrators to optimal control problems is referred to as computational geometric optimal control.
This formulation is based on the concept of the reachability set on a \Poincare section.
This method allows for one to easily determine potential transfer opportunities by finding set intersections on a lower dimensional space and greatly reduces the design process.
The addition of continuous low thrust propulsion extends the control free design process presented in Reference~\citenum{koon2000} and allows for a greater range of potential transfers with a reduced time of flight.

Reachability theory provides a framework to evaluate control capability and safety.  
The reachable set contains all possible trajectories that are achievable over a fixed time horizon from a defined initial condition, subject to the operational constraints of the system.
Reachability theory has been applied to aerospace systems such as collision avoidance, safety planning, and performance characterization.
The theory formally supporting reachability has been extensively developed and is directly derivable from optimal control theory~\cite{varaiya2000,lygeros2002,lygeros2004}.
Computation of the reachable set for a system involves solving the Hamilton-Jacobi partial differential equation or satisfying a dynamic programming principle.
Analytical computation of reachable sets is an ongoing problem and is only possible for certain classes of systems.
Typically, numerical methods are used to generate approximations of the reachability set, but are generally limited by the dimensionality of the problem.
 
Reachability theory has recently been applied to space systems~\cite{holzinger2009,komendera2012a}.
Computation of reachable sets is critical to space situational awareness, rendezvous and proximity operations, and orbit determination operations.
Specifically maintaining accurate estimates of a spacecraft state over extended periods is not trivial.
The challenge is increased for multiple spacecraft operating in close proximity or when there are long periods of time between measurements.
Coupling the ability for continuous low-thrust propulsion between measurements increases the measurement association complexity.
Computing the reachability set given estimated states and control authorities allows one to better correlate subsequent measurements or determine sensor pointing regions in the event of a lost spacecraft. 

The cost function is defined as
\begin{equation}
	J = -\frac{1}{2} \left( \bar{x}(N) - \bar{x}_{n}(N)\right)^T 
	\left[
	\begin{array}{cccc}
		1 & 0& 0& 0 \\
		 0& 0& 0& 0\\
		 0 & 0 & 1 &0\\
		 0 & 0& 0& 0
	\end{array}
	\right]
	\left( \bar{x}(N) - \bar{x}_{n}(N)\right) \, .
	\label{eq:cost}
\end{equation}
The term \( \bar{x}_n(N) \) is the final state of a control-free trajectory while the term \( \bar{x}(N) \) is the final state under the influence of the control input.
In this fashion, the aim is to maximize the distance of the final state from that of the control-free trajectory. 
A chosen \Poincare section is defined through the use of appropriate terminal constraints given by
\begin{subequations}
\begin{align}
    m_1 &= 0 = \frac{y(N) - L_{1y}}{x(N) - L_{1x}} - \tan{\alpha_d} \, , \\ 
    m_2&= 0 = \frac{\dot{x}(N) - \dot{x_n}(N) }{x(N) -x_n(N) } - \tan{\theta_d} \, , \\
	 0 &\geq\bar{u}^T \bar{u} - u_{max}^2 \, ,
\end{align}
    \label{eq:constraints}
\end{subequations}
where the angles \( \alpha_d\) and \( \theta_d\) define the \Poincare section and a specific direction upon the section, respectively. 
The goal is to determine the control input \( \bar{u}_k\) such that the cost function~\cref{eq:cost} is minimized subject to the state equations of motion~\cref{eq:discrete_eoms} and constraints~\cref{eq:constraints}.
Application of the Euler-Lagrange equations allows us to derive the necessary conditions for optimality~\cite{bryson1975}.
The discrete variational integrator in~\cref{eq:discrete_eoms} is used rather than the continuous time counterpart. 
This results in a discrete optimal control problem and the discrete necessary conditions are given as
\begin{subequations}\label{eq:necc_cond}
\begin{align}
	\lambda_{k+1}^T &= \lambda_k^T  \deriv{f_k}{\bar{x}_k}^{-1} \, , \\
	0 &=  \deriv{H_k}{\bar{u}_k} \, ,\\
	0 &= \deriv{\phi}{\bar{x}_k}^T + \deriv{m}{\bar{x}_k}^T\beta  - \lambda^T(N) \, ,  
\end{align}
\end{subequations}
where the Hamiltonian \(H\) is defined as
\begin{equation}
	H_k = \lambda_k^T f(\bar{x}_k, \bar{u}_k) \, ,
	\label{eq:hamiltonian_opt}
\end{equation}
and \(\lambda \in \R^{4 \times 1}\) is the costate and \(\beta \in \R^{2 \times 1} \) are the additional Lagrange multipliers associated with the terminal constraints in~\cref{eq:constraints}.
This indirect optimal control formulation leads to a two point boundary value problem with split boundary conditions. 
By sweeping the angle \( \theta_d \) one can approximate the reachable set on the \Poincare section subject to the bounded control input. 

The costate equation of motion requires the Jacobian of~\cref{eq:discrete_eoms} and is given by
\begin{align}\label{eq:costate_eom}
	\lambda_{k+1}^T = \lambda_k^T
	\begin{bmatrix} 
		f_{1_x} & f_{1_y} & f_{1_{\dot{x}}} & f_{1_{\dot{y}}} \\
		f_{2_x} & f_{2_y} & f_{2_{\dot{x}}} & f_{2_{\dot{y}}} \\
		f_{3_x} & f_{3_y} & f_{3_{\dot{x}}} & f_{3_{\dot{y}}} \\
		f_{4_x} & f_{4_y} & f_{4_{\dot{x}}} & f_{4_{\dot{y}}}
	\end{bmatrix} ^ {-1} \, .
\end{align}
The derivation of~\cref{eq:costate_eom} is given in Appendix A.
In addition, the computation of~\cref{eq:costate_eom} requires inversion of the Jacobian matrix.
This is a computationally expensive operation that is prone to numerical error and instability.
A method is presented in Appendix B to avoid this inversion and determine an explicit update map \( \lambda_k \to \lambda_{k+1} \).

The optimal control formulation presented in this section results in a two point boundary value problem (TPBVP). 
There exist many methods to solve TPBVPs such as gradient, quasilinearization, and shooting methods~\cite{bryson1975}.
In this work, a multiple shooting method is implemented.
Shooting methods are common in astrodynamic trajectory design problems and relatively simple to implement.
In the shooting method, initial conditions are varied such that a terminal constraint is satisfied, similar to the way an archer modifies the bow in order to more accurately strike a target. 
Rather than numerical integration over the entire time interval, multiple shooting segments the interval into several smaller sub-arcs.
Additional interior constraints are used to ensure a continuous solution.
This has the benefit of decreasing the numerical sensitivity of the final states to changes in the initial conditions along each sub-arc.
\section{Numerical Example}\label{sec:simulation}
A numerical simulation is presented to demonstrate the transfer procedure.
The goal is to design a transfer trajectory from a planar periodic orbit about the \( L_1\) Lagrange point to a bounded orbit in the vicinity of the Moon.
The target region is created by choosing an initial condition of \( x_0 = \begin{bmatrix}1.05 & 0 & 0 & 0.35 \end{bmatrix} \) with \( \mu = 0.0125 \).
The initial state is propagated over a period of \( t = \num{20} \) in non-dimensional units which corresponds to approximately \num{1.5} years.
\begin{figure}[htbp]
   \centering
   \includegraphics[width=0.5\textwidth]{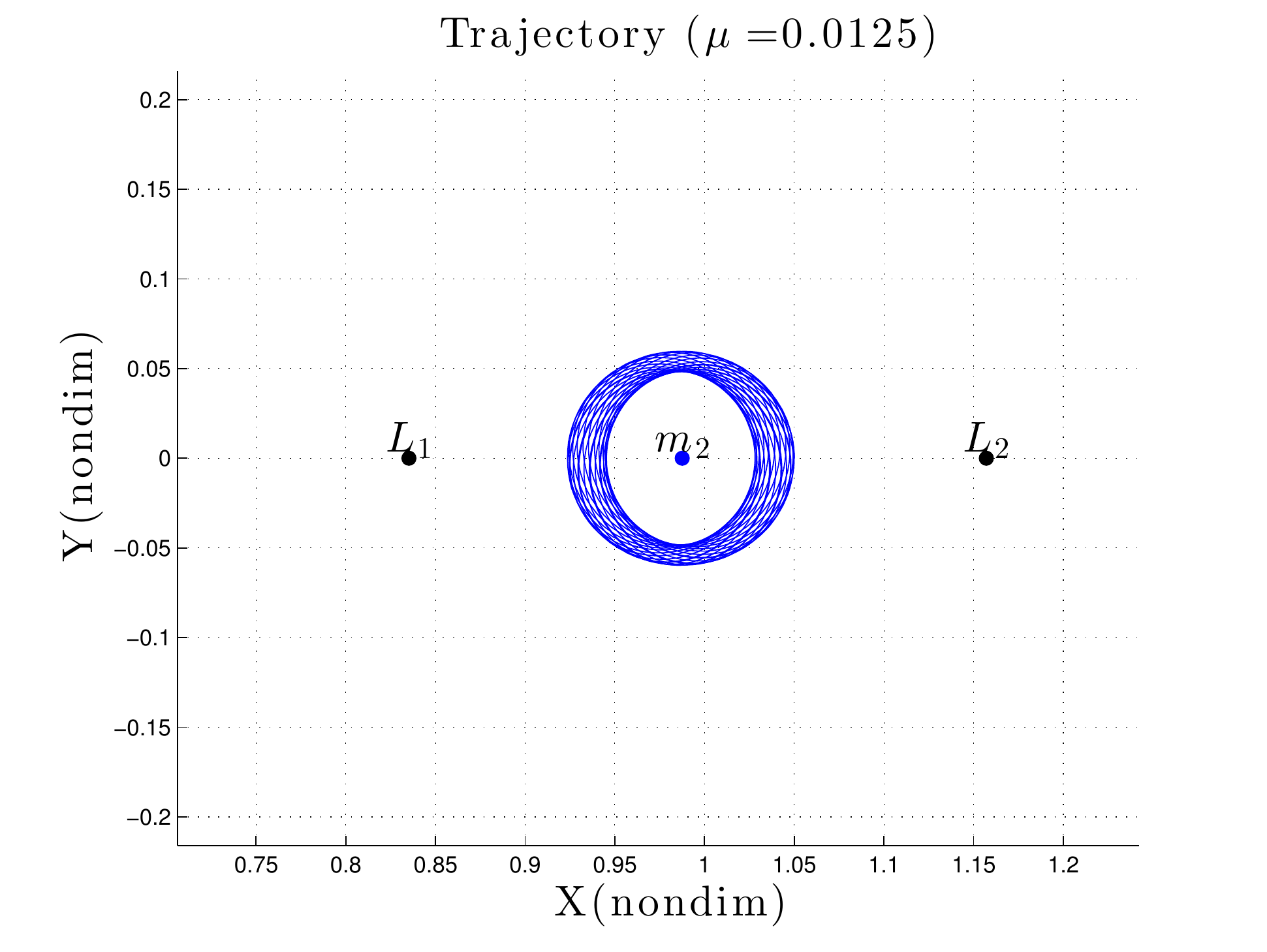} 
   \caption{Target Orbit Region}
   \label{fig:moon_orbit}
\end{figure}
\Cref{fig:moon_orbit} shows that the desired trajectories remain in the vicinity of the Moon in the rotating reference frame. 
This type of orbit would be useful for a variety of mission scenarios.
For example, a series of communication satellites could be placed in this type of orbit. 
The bounded trajectories of the vehicles and constant line of sight to both the Moon and the Earth would allow for constant communication for future manned missions and potential habitats.

\begin{figure} 
	\centering 
	\begin{subfigure}[htbp]{0.5\textwidth} 
		\includegraphics[width=\textwidth]{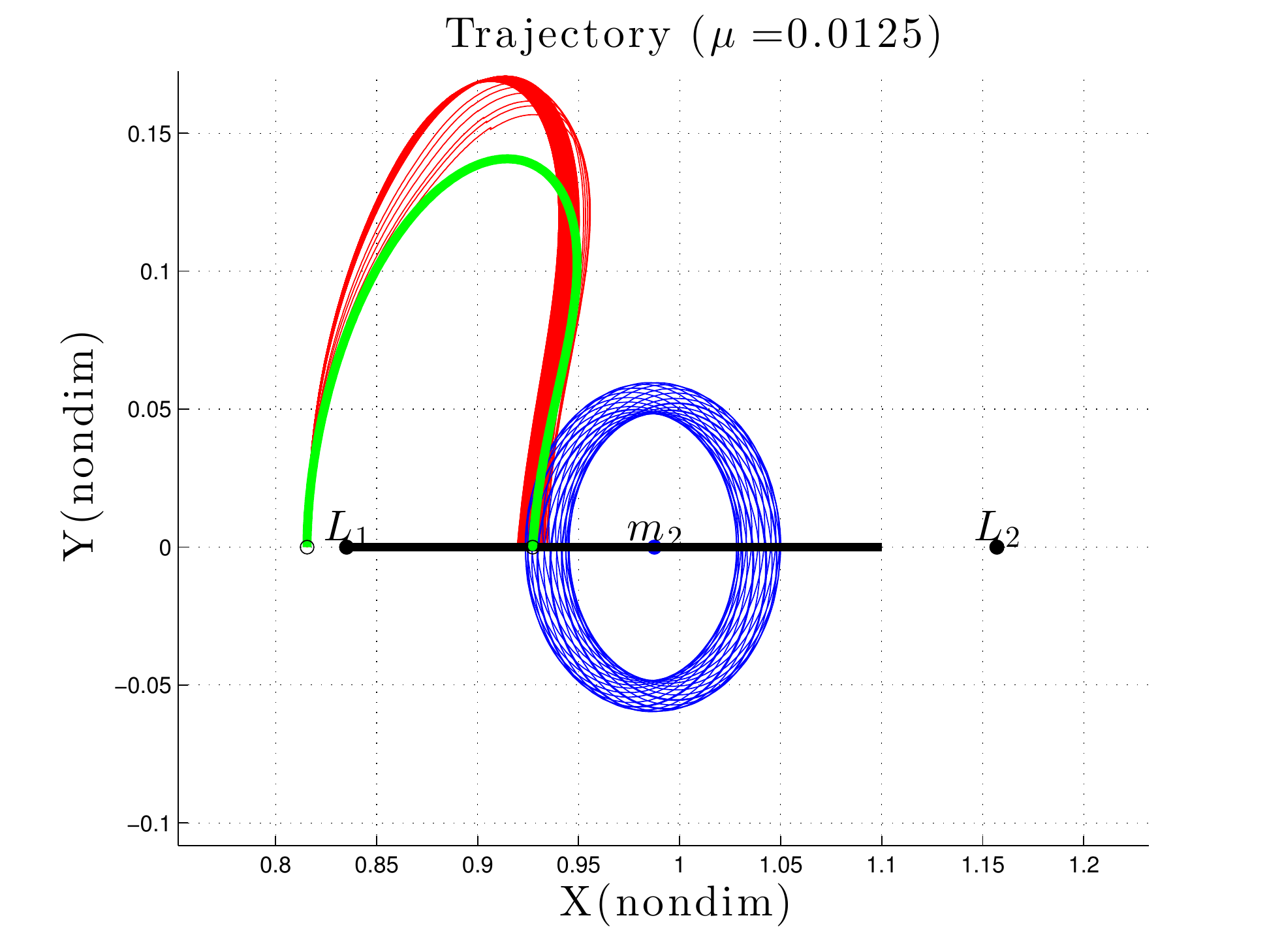} 
		\caption{Controlled Trajectories} \label{fig:reach_trajectory} 
	\end{subfigure}~ 
	\begin{subfigure}[htbp]{0.5\textwidth} 
		\includegraphics[width=\textwidth]{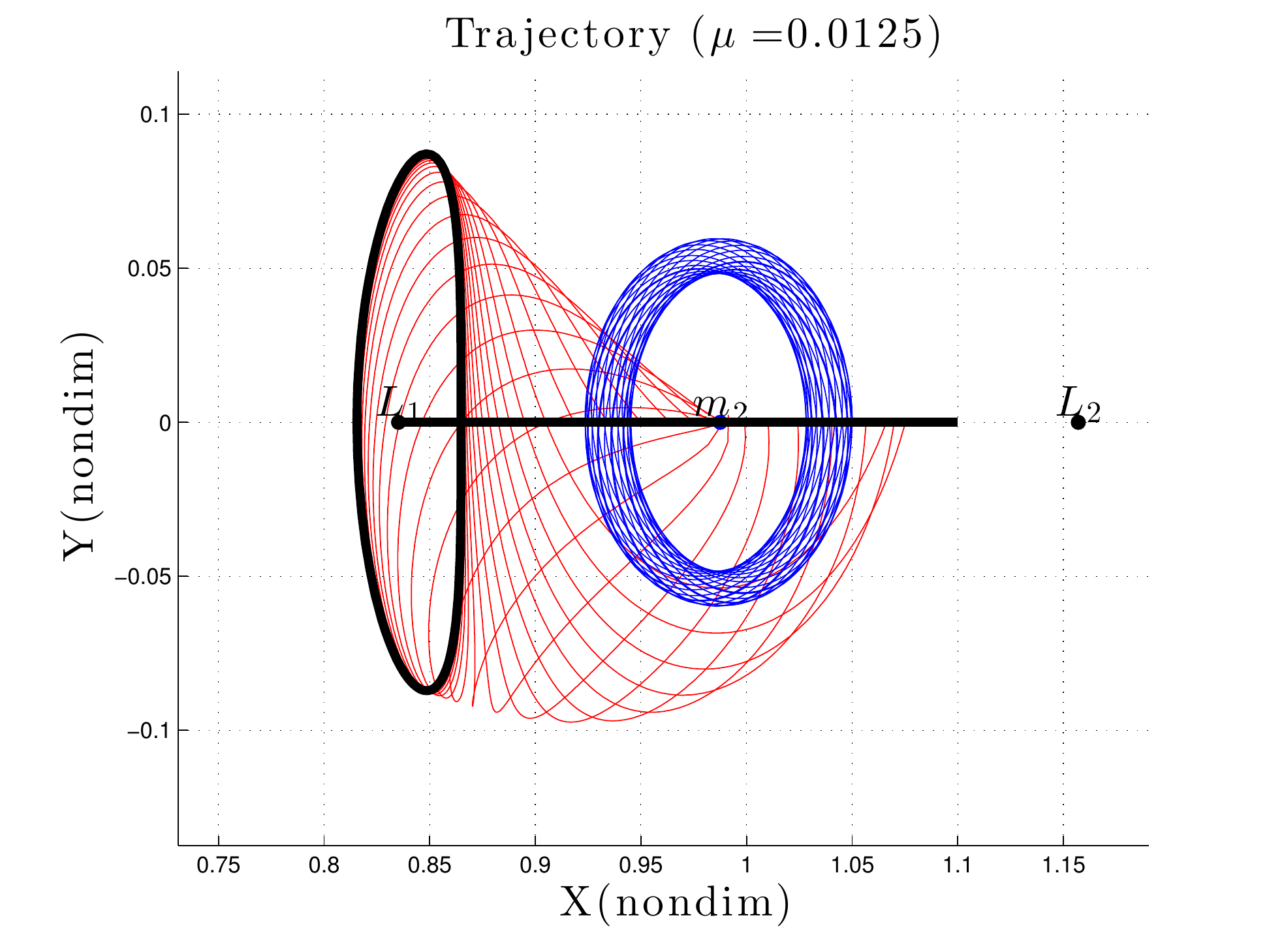} 
		\caption{Invariant Manifold} \label{fig:manifold_trajectory} 
	\end{subfigure} 
	
	\begin{subfigure}[htbp]{0.5\textwidth} 
		\includegraphics[width=\textwidth]{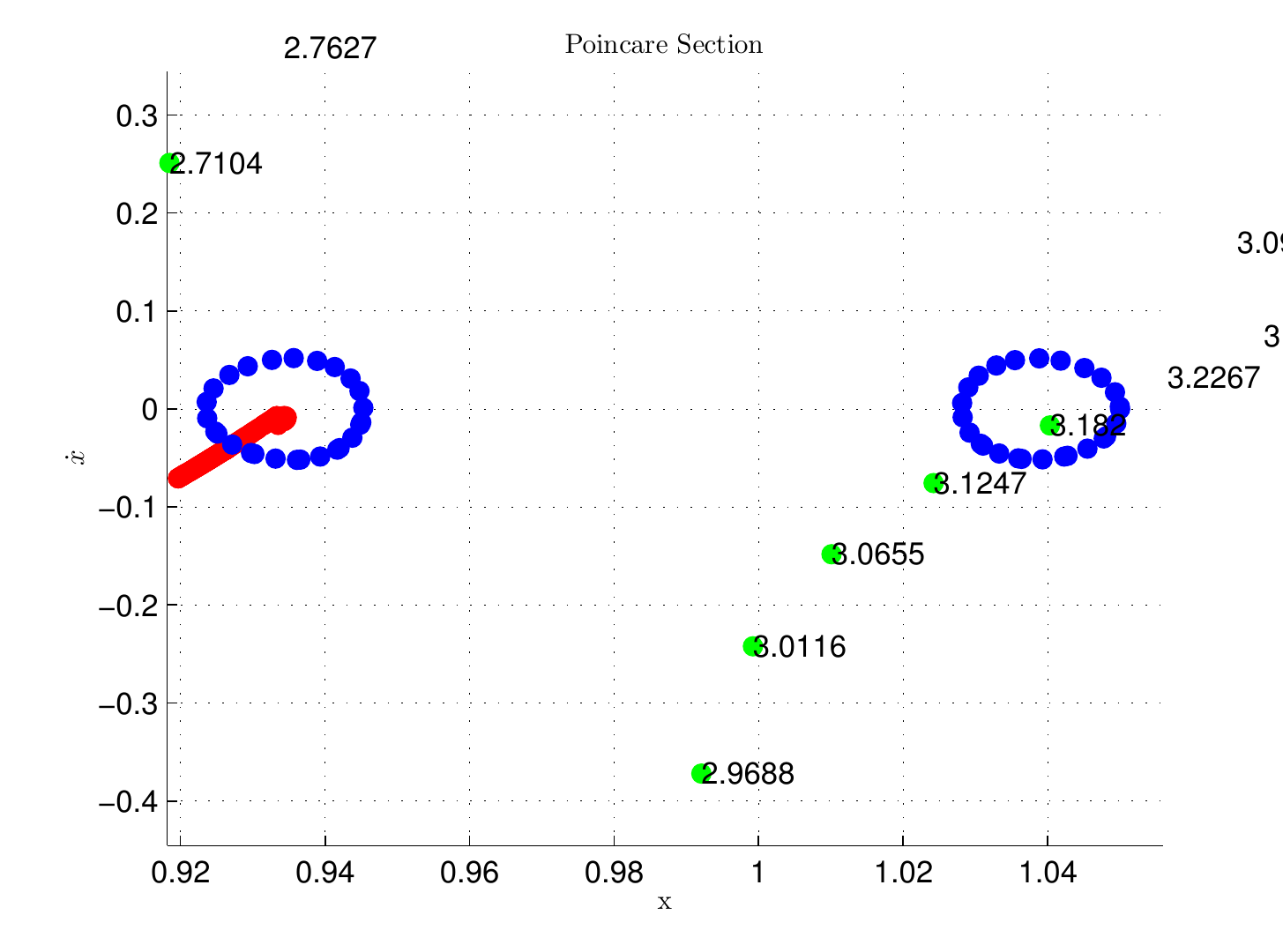} 
		\caption{Poincar\`e Section} \label{fig:poincare_compare} 
	\end{subfigure} 
	\caption{Transfer Trajectory}
	\label{fig:transfer_orbit} 
\end{figure}

As a baseline, the method introduced in Reference~\citenum{koon2000} is implemented.
The method is based on the invariant manifolds associated with the periodic orbits of the three body system.
These invariant manifolds are a set of trajectories that either asymptotically arrive or depart the periodic orbit. 
Linking the invariant manifolds of several periodic orbits, or coupled three body systems, has been used to generate orbital transfers.
\Cref{fig:manifold_trajectory} shows the unstable invariant manifold associated with the \( L_1\) periodic orbit. 
The invariant manifolds are globalized using the eigenvectors of the Monodromy matrix along the periodic orbit.
The eigenvectors serve as a local approximation of the stable and unstable directions. 
Perturbations along these directions serve to approximate the invariant manifolds.

The unstable invariant manifold is numerically propagated to the same \Poincare section defined on the \( \hat{x} \) axis.
The intersection of the invariant manifold and the \Poincare section is denoted by the green markers in~\cref{fig:poincare_compare}.
Only a single branch of the invariant manifold intersects with the ascending region of the target orbit.
There are no intersections of the invariant manifold with the descending region of the target orbit.
The numerical values associated with the green points denote the required time of flight along the invariant manifold in non-dimensional units.
A transfer along the invariant manifold requires on average \( t_f \approx 3.1 \) as compared to \( t_f \approx 1.4 \) for a transfer using low thrust propulsion and the reachable set.
Finally, it should be noted that an additional maneuver would be required for a transfer using the invariant manifold.
The intersection on the \Poincare section only shows that the \( x \text{ and } \dot{x} \) components intersect.
An additional instantaneous \( \Delta V \) would be required to transfer from the invariant manifold to the target.

A periodic orbit is generated about the \( L_1 \) Lagrange point. 
A \Poincare section is chosen to allow for design on a lower dimensional subspace.
The section is defined along the \( \hat{x} \) axis, or defined with \( \alpha = \ang{0} \) and intersects both the initial and target orbits.
From the periodic orbit, a series of optimal trajectories are generated to approximate the reachable set.
The trajectories are generated from a fixed initial state of \( x_0 = \begin{bmatrix}0.8156 & 0 & 0 & 0.1922 \end{bmatrix} \) over a fixed time span of \( t_f = 1.4 \).
By varying the angle \( \theta_d\) in~\cref{eq:constraints}, a different direction along the \Poincare section is maximized. 
The intersection of the optimal trajectories as well as those of the target Moon orbit with the \Poincare section are shown in~\cref{fig:transfer_orbit}.
The optimal trajectories, under the influence of the control input \( \bar{u} \), are shown in~\cref{fig:reach_trajectory}.
Initially, the spacecraft is assumed to lie on the periodic orbit.
As a result, the intersection of this periodic orbit with the \Poincare section are two points corresponding to the two crossing of the orbit.
The use of the continuous low thrust propulsion expands the reachable set to region bounded by the red markers in~\cref{fig:poincare_compare}.
The reachable set is an ellipsoidal region with a major axis aligned along \( \theta \approx \ang{70} \).

The blue points in~\cref{fig:poincare_compare} are the intersections of the target Moon orbit and the \Poincare section.
The two circular regions are the ascending (right) and descending (left) intersections of the target orbit and \Poincare section.
\Cref{fig:poincare_compare} shows that the reachable set and those of the descending target region intersect.
As both regions are discretely approximated a linear interpolation is used to determine the exact intersection state on the \Poincare section.
This intersection generates a partial target state of \( x_t \text{ and } \dot{x}_t \).
Using~\cref{eq:jacobi} and the intersection state the final component \( \dot{y} \) is calculated. 
This results in a complete target state \( \bar{x}_t \) which lies on the reachable set and on the target orbit. 
A final optimal trajectory is generated such that the \( \bar{x}(N) = \bar{x}_t \).
This transfer trajectory is denoted by the green path in~\cref{fig:reach_trajectory}.

%
%
%
%
\section{Conclusions}\label{sec:conclusion}
In this paper, an optimal transfer process which combines concepts of reachability and \Poincare section is used to generate transfer between planar periodic orbits in the three-body problem.
The \Poincare section allows for trajectory design on a lower dimensional phase space and simplifies the process.
The indirect optimal control formulation enables straightforward method of incorporating additional path and control constraints.
However, the use of optimal control techniques leads to open loop trajectories that are not robust to model uncertainties or disturbances.
Lyapunov control theory, which has previously been applied to the two-body problem, is being investigated in the hope of designing closed loop control schemes for this three-body scenario~\cite{chang2002}.
This analysis has also assumed perfect attitude control and the ability to orient thrust in any direction.
The addition of attitude dynamics and realistic pointing constraints would significantly improve the applicability.

\bibliographystyle{AAS_publication}   
\bibliography{library}

\begin{thebibliography}{10}

\bibitem{haque2013}
S.~E. Haque, M.~Keidar, and T.~Lee, ``Low-Thrust Orbital Maneuver Analysis for
  Cubesat Spacecraft with a Micro-Cathode Arc Thruster Subsystem,''  {\em
  Proceedings of the thirty-third international electric propulsion conference,
  Electric Rocket Propulsion Society, Washington DC, USA}, 2013.

\bibitem{koon2000}
W.~S. Koon, M.~W. Lo, J.~E. Marsden, and S.~D. Ross, {\em Dynamical Systems,
  the Three-Body Problem and Space Mission Design}.
\newblock World Scientific, 2000.

\bibitem{ross2006}
S.~D. Ross, ``The interplanetary transport network,''  {\em American
  scientist}, Vol.~94, No.~3, 2006, p.~230.

\bibitem{gomez2001}
G.~G{\'o}mez, W.~Koon, M.~W. Lo, J.~Marsden, J.~Masdemont, and S.~Ross,
  ``Invariant Manifolds, the Spatial Three-Body Problem and Space Mission
  Design,''  {\em AAS/AIAA Astrodynamics Specialist Conference, Quebec City,
  Canada, 2001}, American Astronautical Society, August 2001.

\bibitem{mingotti2011}
G.~Mingotti, F.~Topputo, and F.~Bernelli-Zazzera, ``Earth--Mars transfers with
  ballistic escape and low-thrust capture,''  {\em Celestial Mechanics and
  Dynamical Astronomy}, Vol.~110, No.~2, 2011, pp.~169--188,
  10.1007/s10569-011-9343-5.

\bibitem{grebow2011}
D.~J. Grebow, M.~T. Ozimek, and K.~C. Howell, ``Design of Optimal Low-Thrust
  Lunar Pole-Sitter Missions,''  {\em The Journal of the Astronautical
  Sciences}, Vol.~58, No.~1, 2011, pp.~55--79, 10.1007/BF03321159.

\bibitem{lanczos1970}
C.~Lanczos, {\em The Variational Principles of Mechanics}, Vol.~4.
\newblock Courier Corporation, 1970.

\bibitem{szebehely1967}
V.~Szebehely, ``Theory of orbits. The Restricted Problem of Three Bodies,''
  {\em New York: Academic Press}, Vol.~1, 1967.

\bibitem{marsden2001}
J.~E. Marsden and M.~West, ``Discrete Mechanics and Variational Integrators,''
  {\em Acta Numerica 2001}, Vol.~10, 2001, pp.~357--514.

\bibitem{greenwood1988}
D.~T. Greenwood, {\em Principles of Dynamics}.
\newblock Prentice-Hall Upper Saddle River, NJ, 1988.

\bibitem{conley1968}
C.~C. Conley, ``Low Energy Transit Orbits in the Restricted Three-Body
  Problem,''  {\em SIAM Journal on Applied Mathematics}, Vol.~16, No.~4, 1968,
  pp.~pp. 732--746.

\bibitem{koon1999}
W.~S. Koon, M.~W. Lo, J.~Marsden, and S.~Ross, ``The Genesis Trajectory and
  Heteroclinic Cycles,''  {\em Astrodynamics 1999}, Vol.~103, No.~Part III,
  1999, pp.~2327--2343.

\bibitem{tanaka2011}
K.~Tanaka and J.~Kawaguchi, ``Low-Thrust Transfer between Jovian Moons using
  Manifolds,''  {\em Spaceflight Mechanics}, Vol.~140, No.~Part 3, 2011,
  pp.~1935--1942.

\bibitem{zanzottera2012}
A.~Zanzottera, G.~Mingotti, R.~Castelli, and M.~Dellnitz, ``Intersecting
  invariant manifolds in spatial restricted three-body problems: Design and
  optimization of Earth-to-halo transfers in the Sun--Earth--Moon scenario,''
  {\em Communications in Nonlinear Science and Numerical Simulation}, Vol.~17,
  2 2012, pp.~832--843, http://dx.doi.org/10.1016/j.cnsns.2011.06.032.

\bibitem{campagnola2012}
S.~Campagnola, P.~Skerritt, and R.~Russell, ``Flybys in the planar, circular,
  restricted, three-body problem,''  Vol.~113, No.~3, 2012, pp.~343--368,
  10.1007/s10569-012-9427-x.

\bibitem{ozimek2010a}
M.~T. Ozimek and K.~C. Howell, ``Low-Thrust Transfers in the Earth-Moon System,
  Including Applications to Libration Point Orbits,''  {\em Journal of
  Guidance, Control, and Dynamics}, Vol.~33, 2014/08/16 2010, pp.~533--549,
  10.2514/1.43179.

\bibitem{mingotti2009}
G.~Mingotti, F.~Topputo, and F.~Bernelli-Zazzera, ``Low-energy, low-thrust
  transfers to the Moon,''  {\em Celestial Mechanics and Dynamical Astronomy},
  Vol.~105, No.~1-3, 2009, pp.~61--74.

\bibitem{koon2001}
W.~S. Koon, M.~W. Lo, J.~E. Marsden, and S.~D. Ross, ``Low energy transfer to
  the Moon,''  {\em Dynamics of Natural and Artificial Celestial Bodies},
  pp.~63--73, Springer, 2001.

\bibitem{varaiya2000}
P.~Varaiya, ``Reach set computation using optimal control,''  {\em Verification
  of Digital and Hybrid Systems}, pp.~323--331, Springer, 2000.

\bibitem{lygeros2002}
J.~Lygeros, ``On the Relation of Reachability to Minimum Cost Optimal
  Control,''  {\em Decision and Control, 2002, Proceedings of the 41st IEEE
  Conference on}, Vol.~2, Dec 2002, pp.~1910--1915 vol.2,
  10.1109/CDC.2002.1184805.

\bibitem{lygeros2004}
J.~Lygeros, ``On reachability and minimum cost optimal control,''  {\em
  Automatica}, Vol.~40, 6 2004, pp.~917--927,
  http://dx.doi.org/10.1016/j.automatica.2004.01.012.

\bibitem{holzinger2009}
M.~Holzinger and D.~Scheeres, ``Reachability Analysis Applied to Space
  Situational Awareness,''  {\em Advanced Maui Optical and Space Surveillance
  Technologies Conference}, 2009.

\bibitem{komendera2012a}
E.~E. Komendera, D.~J. Scheeres, and E.~Bradley, ``Intelligent Computation of
  Reachability Sets for Space Missions.,''  {\em IAAI}, 2012.

\bibitem{bryson1975}
A.~E. Bryson and Y.~C. Ho, {\em Applied Optimal Control: Optimization,
  Estimation and Control}.
\newblock CRC Press, 1975.

\bibitem{chang2002}
D.~E. Chang, D.~F. Chichka, and J.~E. Marsden, ``Lyapunov-based Transfer
  between Elliptic Keplerian Orbits,''  {\em Discrete and Continuous Dynamical
  Systems Series B}, Vol.~2, No.~1, 2002, pp.~57--68.

\end{thebibliography}

\appendix
\section*{Appendix A: Costate Equations of Motion}\label{sec:costate_appendix}
The development of the costate equations of motions begins with determining the second order partial derivatives of the gravitational potential. 
Due to the symmetry of partial derivatives only three terms are required and are given by
\begin{align}\label{eq:second_discrete_potential_grad}
	U_{x\xk} &= \parenth{1-\mu} \bracket{\frac{1}{\distonek^3} - \frac{3 \parenth{\xk +\mu}^2}{\distonek^5}} + \mu \bracket{\frac{1}{\disttwok^3} - \frac{3 \parenth{\xk -1 + \mu}^2}{\disttwok^5}} \, ,\\
	U_{y\yk} &= \parenth{1-\mu} \bracket{\frac{1}{\distonek^3} - \frac{3 \yk^2}{\distonek^5}} + \mu \bracket{\frac{1}{\disttwok^3} - \frac{3 \yk^2}{\disttwok^5}} \, ,\\
	U_{x\yk} &= U_{y\xk} =  \frac{-3 \parenth{1-\mu} \parenth{\xk +\mu} \yk}{\distonek^3} - \frac{3\mu\yk\parenth{\xk-1+\mu}}{\disttwok^5} \, .
\end{align}
The gradient of~\cref{eq:xkp} is given as
\begin{subequations}\label{eq:xkp_grad}
\begin{align}
	f_{1_x} &= \frac{1}{1+h^2} \bracket{h^2 + 1 + \frac{h^2}{2} -\frac{h^3}{2} U_{y\xk} - \frac{h^2}{2}U_{x\xk}} \, ,\\
	f_{1_y} &= \frac{1}{1+h^2} \bracket{ \frac{h^3}{2} -\frac{h^3}{2} U_{y\yk} - \frac{h^2}{2}U_{x\yk}} \, ,\\
	f_{1_{\dot{x}}} &= \frac{h}{1+h^2} \, ,\\
	f_{1_{\dot{y}}} &= \frac{h^2}{1+h^2} \, .
\end{align}
\end{subequations}
The gradient of~\cref{eq:ykp} is given as
\begin{subequations}\label{eq:ykp_grad}
\begin{align}
	f_{2_x} &= h -h f_{1_x} - \frac{h^2}{2} U_{y\xk} \, , \\
	f_{2_y} &= -h f_{1_y} + 1 + \frac{h^2}{2} - \frac{h^2}{2} U_{y\yk} \, ,\\
	f_{2_{\dot{x}}} &= -h f_{1_{\dot{x}}} \, ,\\
	f_{2_{\dot{y}}} &= h - h f_{1_{\dot{y}}} \, .
\end{align}
\end{subequations}
The gradients of~\cref{eq:distonekp,eq:disttwokp} are given as
\begin{subequations}\label{eq:distkp_grad}
\begin{align}
	\deriv{\distonekp}{\bar{x}} &= \parenth{\left( \xkp + \mu\right)^2 + \ykp^2}^{-\frac{1}{2}} \bracket{\parenth{\xkp + \mu} f_{1_{\bar{x}}} + \ykp f_{2_{\bar{x}}}} \, ,\\
	\deriv{\disttwokp}{\bar{x}} &= \parenth{\left( \xkp - 1 + \mu\right)^2 + \ykp^2}^{-\frac{1}{2}} \bracket{\parenth{\xkp -1 + \mu} f_{1_{\bar{x}}} + \ykp f_{2_{\bar{x}}}}  \, .
\end{align}
\end{subequations}
The second order partial derivatives of the gravitational potential at \( k+1\) are given as
\begin{subequations}
\begin{align}
	\deriv{U_{x\xkp}}{\bar{x}} &= \parenth{1-\mu}\bracket{\frac{1}{\distonekp^3} f_{1_{\bar{x}}} - \frac{3 \parenth{\xkp +mu}}{\distonekp^4} \deriv{\distonekp}{\bar{x}}} + \mu \bracket{\frac{1}{\disttwokp^3} f_{1_{\bar{x}}} - \frac{-3 \parenth{\xkp -1 + \mu}}{\disttwokp^4} \deriv{\disttwokp}{\bar{x}}} \, ,\\
	\deriv{U_{y\xkp}}{\bar{x}} &= \parenth{1-\mu}\bracket{\frac{1}{\distonekp^3} f_{2_{\bar{x}}} - \frac{3 \ykp}{\distonekp^4} \deriv{\distonekp}{\bar{x}}} + \mu \bracket{\frac{1}{\disttwokp^3} f_{2_{\bar{x}}} - \frac{-3 \ykp}{\disttwokp^4} \deriv{\disttwokp}{\bar{x}}} \, .
\end{align}
\end{subequations}
The gradient of~\cref{eq:xdotkp,eq:ydotkp} are given as
\begin{subequations}\label{eq:xdotkp_grad}
\begin{align}
	f_{3_x} &= 2 f_{2_x} + \frac{h}{2} \parenth{f_{1_x} + 1} - \frac{h}{2} U_{x\xkp} - \frac{h}{2} U_{x\xk} \, ,\\
	f_{3_y} &= -2 + 2 f_{2_y} + \frac{h}{2} f_{1_y} - \frac{h}{2} U_{x\ykp} -\frac{h}{2} U_{x\yk} \, ,\\
	f_{3_{\dot{x}}} &= 1 + 2 f_{2_{\dot{x}}} + \frac{h}{2} f_{1_{\dot{x}}} - \frac{h}{2} U_{x\xdotkp}\, ,\\
	f_{3_{\dot{y}}} &= 2 f_{2_{\dot{y}}}\, ,
\end{align}
\end{subequations}
\begin{subequations}\label{eq:ydotkp_grad}
\begin{align}
	f_{4_x} &= 2 - 2 f_{1_x} + \frac{h}{2} f_{2_x}  - \frac{h}{2} U_{y\xkp} - \frac{h}{2} U_{y\xk} \, ,\\
	f_{4_y} &= -2f_{1_y}  - \frac{h}{2} \parenth{f_{2_y} + 1} - \frac{h}{2} U_{y\ykp} -\frac{h}{2} U_{y\yk} \, ,\\
	f_{4_{\dot{x}}} &= - 2 f_{1_{\dot{x}}} +  \frac{h}{2} f_{2_{\dot{x}}} - \frac{h}{2} U_{y\xdotkp}\, ,\\
	f_{4_{\dot{y}}} &= 1 - 2 f_{1_{\dot{y}}} + \frac{h}{2} f_{2_{\dot{y}}} - \frac{h}{2} U_{y\ydotkp}\, .
\end{align}
\end{subequations}
These gradient equations are in a cascade type structure.
\Cref{eq:ydotkp_grad,eq:xdotkp_grad} are functions of~\cref{eq:xdotkp,eq:ydotkp}.
As a result, the accuracy of the Jacobian will tend to decrease as the first order approximation errors accumulate.

\section*{Appendix B: Gauss Jordan Elimination}\label{sec:costate_gauss_jordan}
The costate equations of motion are given by~\cref{eq:costate_eom} and repeated here as
\begin{align}\label{eq:costate_eom_transpose}
	\begin{bmatrix} 
		\fonex & \ftwox & \fthreex & \ffourx \\
		\foney & \ftwoy & \fthreey & \ffoury \\
		\fonexd & \ftwoxd & \fthreexd & \ffourxd \\
		\foneyd & \ftwoyd & \fthreeyd & \ffouryd
	\end{bmatrix}
	\begin{bmatrix} \lambda_{\xkp} \\ \lambda_{\ykp} \\ \lambda_{\xdotkp} \\ \lambda_{\ydotkp} \end{bmatrix}
	=
	\begin{bmatrix} \lambda_{\xk} \\ \lambda_{\yk} \\ \lambda_{\xdotk} \\ \lambda_{\ydotk} \end{bmatrix} \, .
\end{align}
To determine the discrete update map \( \lambda_k \to \lambda_{k+1}\) the inverse of the Jacobian matrix is required.
In order to avoid the need of an explicit inversion a Gauss Jordan method is implemented.
To begin, several terms are defined which are required to carry out the row operations and are defined as
\begin{subequations}\label{eq:ref_scale}
\begin{align}
	a &= -\frac{\foney}{\fonex} \, ,\\
	b &= -\frac{\fonexd}{\fonex} \, , \\
	c &= -\frac{\foneyd}{\fonex} \, ,\\
	e &= -\frac{\ftwoxd + b \ftwox}{\ftwoy + a \ftwox} \, ,\\
	f &= -\frac{\ftwoyd + c \ftwox}{\ftwoy + a \ftwox} \, ,\\
	g &= -\frac{\fthreeyd + c \fthreex + f \parenth{\fthreey + a \fthreex}}{\fthreexd + b \fthreex + e \parenth{\fthreey + a \fthreex}}\, .
\end{align}
\end{subequations}
\Cref{eq:costate_eom_transpose} is transformed to row echelon form using elementary row operations and is defined as
\begin{align}\label{eq:costate_ref}
	\begin{bmatrix} 
		\alpha_{11} & \alpha_{12} & \alpha_{13} & \alpha_{14} \\
		0 & \alpha_{22} & \alpha_{23} & \alpha_{24} \\
		0 & 0 & \alpha_{33} & \alpha_{34} \\
		0 & 0 & 0 & \alpha_{44}
	\end{bmatrix}
	\begin{bmatrix} \lambda_{\xkp} \\ \lambda_{\ykp} \\ \lambda_{\xdotkp} \\ \lambda_{\ydotkp} \end{bmatrix}
	=
	\begin{bmatrix} \beta_1 \\ \beta_2 \\ \beta_3 \\ \beta_4 \end{bmatrix} \, ,
\end{align}
where the terms \( \alpha_{ij} \) and \( \beta_{i} \) are defined as follows
\begin{subequations}
\begin{align}
	\alpha_{11} &= \fonex \, ,\\
	\alpha_{12} &= \ftwox \, ,\\
	\alpha_{13} &= \fthreex \, ,\\
	\alpha_{14} &= \ffourx \, ,\\
	\alpha_{22} &= \ftwoy + a \ftwox \, ,\\
	\alpha_{23} &= \fthreey + a\fthreex \, ,\\
	\alpha_{24} &= \ffoury + a \ffourx \, ,\\
	\alpha_{33} &= \fthreexd + b \fthreex + e \parenth{\fthreey + a \fthreex}\, ,\\
	\alpha_{34} &= \ffourxd + b \ffourx + e \parenth{\ffoury + a \ffourx} \, ,\\
	\alpha_{44} &= \ffouryd + c \ffourx + f \parenth{\ffoury + a \ffourx} + g \parenth{\ffourxd + b \ffourx + e \parenth{\ffoury + a \ffourx}} \, ,\\
	\beta_1 &= \lambda_{\xk} \, ,\\
	\beta_2 &= \lambda_{\yk} + a \lambda_{\xk} \, ,\\
	\beta_3 &= \lambda_{\xdotk} + b \lambda_{\xk} + e \parenth{\lambda_{\yk} + a \lambda_{\xk}} \, ,\\
	\beta_4 &= \lambda_{\ydotk} + c \lambda_{\xk} + f \parenth{\lambda_{\yk}+a \lambda_{\xk}} + g \parenth{\lambda_{\xdotk} + b \lambda_{\xk} + e \parenth{\lambda_{\yk} + a \lambda_{\xk}}} \, .
\end{align}
\end{subequations}
Finally, backsubstituion is used to determine explicit equations for the discrete update map \( \lambda_k \to \lambda_{k+1} \) which is defined as
\begin{subequations}\label{eq:costate_update}
\begin{align}
	\lambda_{\ydotkp} &= \frac{\beta_4}{\alpha_{44}} \, ,\\
	\lambda_{\xdotkp} &= \frac{\beta_3}{\alpha_{33}} - \frac{\alpha_{34}}{\alpha_{33}} \lambda_{\ydotkp} \, ,\\
	\lambda_{\ykp} &= \frac{\beta_2}{\alpha_22} - \frac{\alpha_{33}}{\alpha_{22}}\lambda_{\xdotkp} - \frac{\alpha_{24}}{\alpha_{22}} \lambda_{\ydotkp} \, ,\\
	\lambda_{\xkp} &= \frac{\beta_1}{\alpha_{11}} - \frac{\alpha_{12}}{\alpha_{11}} \lambda_{\ykp} - \frac{\alpha_{13}}{\alpha_{11}} \lambda_{\xdotkp} - \frac{\alpha_{14}}{\alpha_{11}} \lambda_{\ydotkp} \, .
\end{align}
\end{subequations}

\end{document}